\documentclass[onefignum,onetabnum]{siamart190516}
\usepackage[utf8]{inputenc}
\usepackage{amssymb}
\usepackage{float}
\usepackage{tikz,tikz-cd}
\usepackage{amsmath}
\usepackage{algorithmic}
\usetikzlibrary{cd}
\usepackage{xcolor}
\pagecolor{white}

\usepackage[normalem]{ulem} 
\newcommand{\strikeout}[1]{{\color{red}{\expandafter\sout\expandafter{#1}}}}

\def\bv{\mathbf{v}}
\def\bn{\mathbf{n}}
\def\bt{\mathbf{t}}

\DeclareMathOperator{\ima}{im}

\title{Combinatorial   and Hodge Laplacians: \\ Similarity and Difference}
\author{Emily Ribando-Gros, Rui Wang, Jiahui Chen, Yiying Tong, Guo-Wei Wei}

\begin{document}

\maketitle

\begin{abstract}
As key subjects in spectral geometry and combinatorial graph theory respectively, the (continuous) Hodge Laplacian and the combinatorial  Laplacian share similarities in revealing the topological dimension and geometric shape of data and in their realization of diffusion and minimization of harmonic measures.  It is believed that they also both associate with vector calculus, through the gradient, curl, and divergence, as argued  in the popular usage of ``Hodge Laplacians on graphs'' in the literature. Nevertheless, these Laplacians are intrinsically different in their domains of definitions and applicability to specific data formats, hindering any in-depth comparison of the two approaches.
For example, the spectral decomposition of a vector field  on a simple point cloud  using combinatorial Laplacians defined on some commonly used simplicial complexes
does not give rise to the same curl-free and divergence-free components as one would obtain from the  spectral decomposition of a vector field  using either the continuous Hodge Laplacians defined on differential forms in manifolds or the discretized   Hodge Laplacians defined on a point cloud with boundary in the Eulerian representation or on a regular mesh in the Eulerian representation.
To facilitate the comparison and bridge the gap between the combinatorial   Laplacian and Hodge Laplacian for the discretization of continuous manifolds with boundary, we further introduce Boundary-Induced Graph (BIG) Laplacians using tools from Discrete Exterior Calculus (DEC).  BIG Laplacians are defined on discrete domains with appropriate boundary conditions to characterize the topology and shape of data.  The similarities and differences of the combinatorial Laplacian, BIG Laplacian, and Hodge Laplacian are then examined. Through an Eulerian representation of 3D domains as level-set functions on regular grids, we show experimentally the conditions for the convergence of BIG Laplacian eigenvalues to those of the Hodge Laplacian for elementary shapes.

\end{abstract}

\section{Introduction}

The Laplacian operator $\Delta$ is ubiquitous, proving a useful formalism in many scientific fields, including spectral geometry, differential geometry, spectral graph theory, and algebraic topology. Simply stated, it computes how a function deviates from average values in a local neighborhood. Classically, this operator is used to study problems such as steady state, heat diffusion, and wave propagation but has since found many other applications. What are particularly interesting are  topological Laplacian operators that are able to describe the algebraic topology structures in their kernel space \cite{eckmann45, mohar1991laplacian,merris1995survey,ghrist2020cellular}.  For example, persistent combinatorial Laplacian    \cite{wanggraph} is related to persistent homology and Topological Data Analysis (TDA) \cite{zomorodian2005computing,edelsbrunner2008persistent}, a hot topic in data science.

The de Rham-Hodge theory provides both topological and geometric characteristics of manifolds via tools from differential geometry. The relationship between the de Rham-Hodge Laplacian (or Hodge Laplacian) and de Rham-Hodge theory is given by an isomorphism between the operator's kernel,  harmonic forms, and the cohomology classes of manifolds. The dimension of the kernel of the $k$-Laplacian, $k\in\mathbb{Z}$, also known as the Betti number $\beta_k$, counts the number of $k$-dimensional holes. Specifically, $\beta_0$ is the number of connected components, $\beta_1$ shows the number of holes, and $\beta_2$ represents the number of cavities. The first nonzero eigenvalue, known as the Fiedler value, or principal eigenvalue, in graph theory describes connectedness.
The continuous Hodge Laplacian is a vital tool in spectral geometry. In 1966, Mark Kac  \cite{Kac}, once asked if one can determine the shape of a drum from its sound frequencies. More specifically, the question asks if eigenvalues of the Laplacian can uniquely determine geometric properties. While in general it may not be possible to determine the shape of the drum given the spectra due to isospectral manifolds, this question still stimulates various developments in  algebraic topology and differential geometry. Recently, speaking to the importance of higher frequency spectral values, in \cite{specshape} the connection between the Laplacian spectra and a shape is learned and used to generate an unknown shape from only its Laplacian spectrum. The newly developed evolutionary de Rham-Hodge theory \cite{chenevol} may have a high potential to break the degeneracy of the isometry class of a manifold by generating a family of submanifolds via a filtration.

The field of spectral graph theory also has a rich history distinct from de Rham-Hodge theory \cite{eckmann45,chung1997spectral,ren2021discrete}. With the wide use of applications of graphs in the 1970s and 1980s,   combinatorial Laplacians,  also known as the graph Laplacian when speaking of the lowest order combinatorial Laplacian, gained popularity \cite{eckmann45,goldberg2002combinatorial,horak2013spectra}. Combinatorial Laplacians and their spectra are proved valuable in  the study of molecular stability \cite{chemistry},  electrical networks \cite{networks}, biomolecule analysis~\cite{wanggraph,meng2021persistent,chen2022,proteineng}, neuroscience~\cite{Lee2019HarmonicHA}, ranking~\cite{rank}, deep learning~\cite{dl}, signal processing~\cite{roddenberry2022signal},and many others. Note that some authors used the term discrete/combinatorial Hodge Laplacians to refer to combinatorial Laplacians that are purely based on cell complex connectivity~\cite{grady,schaubSC,roddenberry2022signal} since these operators also admit their own Hodge decompositions. We will use the shorter term combinatorial Laplacians to easily distinguish them from the discretized and the continuous Hodge Laplacians, which we refer to as Hodge Laplacians in the rest of the paper.

The analogy between   combinatorial Laplacians and continuous Hodge Laplacians was demonstrated in their description of the topology and shape of data \cite{lim15}.  In these cases, the underlying
algebraic topology structure is  similar to that of Hodge Laplacians and the actions of combinatorial  Laplacians were interpreted  by
those of corresponding Hodge Laplacians in differential geometry   \cite{lim15}.
However, combinatorial Laplacians are inherently defined on discrete data via simplicial complexes or their generalizations,  while (continuous) Hodge Laplacians are defined on continuous manifolds  via differential forms.
Additionally, in many applications,  Hodge Laplacians deal with manifolds with boundary, which may be needed for physically-based systems \cite{zhao2020rham}. Moreover,
there is no intrinsic notion of domain boundary for combinatorial Laplacians.  Finally, Hodge Laplacians facilitate the Hodge decomposition of a vector field into harmonic, curl-free, and divergence-free components in the traditional differential calculus sense, which has extensive applications in continuum mechanics, including fluid mechanics and solid mechanics. In contrast, a combinatorial Laplacian decomposition, while computable, does not offer these physically relevant components since it, in general, lacks an approximation to the metric tensor fields for measuring vectors, but merely reflects the connections between simplices of different dimensions. For these reasons, combinatorial Laplacians should not be regarded as interchangeable with  Hodge Laplacians.

\begin{figure}[h]
	\label{fig:venn}
	\includegraphics[width=\textwidth]{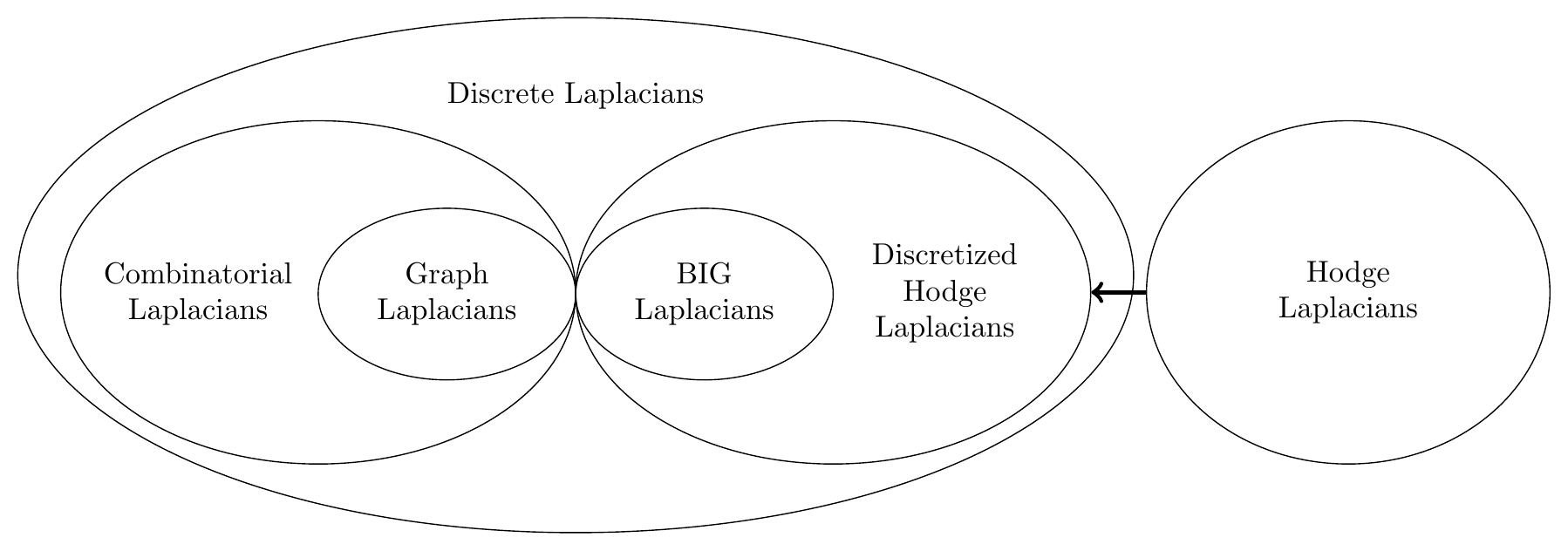}
	\caption{The relationships among different types of topological Laplacians.}
\end{figure}

While  Hodge Laplacians are defined on continuous differential manifolds,  practical computations require discretized  Hodge Laplacians, which are created via
 Discrete Exterior Calculus (DEC \cite{dec06}) or Finite Element Exterior Calculus (FEEC, \cite{fem_ext}), defined on meshes or point clouds with boundary. Even thus, discretized Hodge Laplacians are not compatible to combinatorial Laplacians because they converge to  the original Hodge Laplacians as the largest mesh size approaches to zero.
To put   combinatorial Laplacians and  Hodge Laplacians on an equal footing for the discretization of continuous manifolds with boundary, we introduce   new Boundary-Induced Graph (BIG) Laplacians to facilitate the comparison and contrast of the two concepts.  BIG Laplacians has the discrete nature of combinatorial Laplacians while fully preserving  Hodge Laplacians' capability to perform differential calculus through enforcing appropriate boundary conditions.
The relationships of  various topological Laplacians is illustrated in Fig. \ref{fig:venn}.
Specifically, we discuss Dirichlet and Neumann conditions. A regular grid representation of surfaces also ensures that   BIG Laplacians can numerically deliver results similar to those of  continuous and discretized Hodge Laplacians as shown by our test cases. We further discuss when the eigenvalues of  clique-based combinatorial Laplacians perform unfavorably in this context as opposed to  BIG Laplacians, which are compared to  Hodge Laplacians. Though all operators discussed may be used interchangeably in the simple case of computing Betti numbers, we point out some key differences when seeking higher frequency spectral information, depending on source model and model representation, and when requiring certain boundary conditions. As an independent formulation, beyond illustrating the relation of   combinatorial Laplacians to   Hodge Laplacians,   BIG Laplacians may find its versatile applications due to its numerical benefits.

\subsection{Overview}

 First, in Section \ref{sec:background}, we outline the theory for multidimensional Laplacian operators defined using differential forms. Then, Section \ref{sec:contin} contains a synopsis of the continuous spectral de Rham-Hodge theory with Section \ref{boundarycond} describing the differential form of   Laplacians on bounded domains. Sections \ref{sec:discrete} and \ref{sec:cliquegraph} cover the representations of discrete surfaces, the discretization of the mentioned continuous theories using  DEC, and the usual generic combinatorial Laplacians in this context. After covering the necessary background, we present a novel implementation of BIG  Laplacians for 3D volumes to be comparable to the discretized Hodge Laplacian in Section \ref{sec:big}. Then through simple examples and in-depth numerical calculations in Sections \ref{sec:diff} and \ref{sec:sim}, we explore the different  discretized Hodge Laplacians defined through two standard representations, regular grid Eulerian for volumetric data
and tetrahedral Lagrangian  for volumetric data-induced point cloud with boundary, as well as differing boundary conditions, Dirichlet and Neumann. To summarize, we
\begin{itemize}
	\item introduce   BIG Laplacians alongside   Hodge Laplacians for manifolds with boundary to facilitate the comparison,
	\item illustrate the differences among   usual  combinatorial Laplacians,  new BIG Laplacians as a special case of discretized Hodge Laplacians, and   Hodge Laplacians, and
	\item show the numerical similarities between   BIG Laplacians and Hodge Laplacians.
\end{itemize}

\section{Background of topological Laplacians}\label{sec:background}

In this section, we briefly review a few topological Laplacians.
We discuss the continuous theory on Hodge  Laplacians, boundary conditions for manifolds with boundary, and how  Hodge  Laplacians are commonly discretized on differential forms  through the use of Hodge stars. Then, we recap combinatorial  Laplacians, for which the classical graph Laplacian is a special case.

\subsection{Continuous de Rham-Hodge theory}\label{sec:contin}

Most important physical equations expressed in multivariable calculus can be unified using methods from differential forms and exterior calculus. These methods naturally lead to the de Rham-Hodge theory, in which differential geometry, algebraic topology, partial differential equation, and their relationship to spectral graph theory are revealed for analysis. Fundamentally, de Rham-Hodge theory states that the harmonic part of the spectrum of the Laplacian operator can be used for unique presentations of topological structures, more precisely, the cohomology groups of manifolds. Here a brief overview of the necessary components of continuous de Rham-Hodge theory is provided to define the de Rham-Hodge and combinatorial Laplacians and obtain useful information from their spectral analysis. 

Given a manifold $M$, a differential $k$-form $\omega_k$ is a quantity that can be integrated over a $k$-dimensional domain \footnote{We will try to stick to the notation where the dimension is a subscript except for topological spaces or other special groups, when the dimension will be a superscript. Occasionally the dimension may be left out if clear from the context.}, i.e., $\omega_k$ is an element of $\Omega^k(M)$, the space of antisymmetric rank-$(0,k)$ tensors. For instance, an element of $\Omega^k(\mathbb{R}^n)$ can be seen as a smooth field of linear combinations of antisymmetric $k$-th order tensor products ($\wedge$) of $dx^1,dx^2,\dots,$ and $ dx^n$. In $\mathbb{R}^3$, we may regard 0-forms as scalar fields $f(p)$ that can be evaluated at a point; 1-forms, as linear combinations of $dx^1,$ $dx^2,$ and $dx^3$, can be regarded as a vector field, which can be integrated along arbitrary curves through line integral; 2-forms, as linear combinations of $dx^2\wedge dx^3, dx^3\wedge dx^1,$ and $dx^1\wedge dx^2$, can also be regarded as a vector field, which can be integrated over surface patches as fluxes; 3-forms, as linear combinations of $dx^1\wedge dx^2\wedge dx^3$, can be regarded as a scalar field, which can be integrated over a volume as a density function. In 3D, the wedge products $\wedge$ between 1-forms can be seen as a cross product, and those between 2-forms and 1-forms can be seen as a dot product.

While the following theory holds for surfaces, or more generally Riemannian manifolds, we will mainly consider volumes in order to introduce the concept of boundary conditions in practical problems and increase the applicability of this study. This means that only 0-, 1-, 2-, and 3-forms are considered. From this point, applying certain combinations of exterior calculus operators, such as the exterior derivative and Hodge star, to differential forms unifies the gradient ($\nabla$), divergence ($\nabla\cdot$), curl ($\nabla\times$), and thus the Laplacian ($\Delta, \nabla\cdot\nabla, \nabla^2$) operators in vector field analysis.

In differential geometry, the exterior derivative $d$, also known as the differential operator, explains how quickly a $k$-form changes in every $(k\!+\!1)$-dimensional direction. In fact, $d$ acts as the known gradient, curl, and divergence when applied to 0-, 1-, and 2-forms, respectively.
Furthermore, $d_k:\Omega^k(M)\rightarrow\Omega^{k+1}(M)$ has the essential property that $d_{k+1}\circ d_k=0$. This property provides the links to the topology of the Morse theory on scalar fields, the index theory of singularities of vector fields, and the integrabilities (Hodge decomposition) of vector fields.
Given a metric, the $L_2$ inner product on the space of $k$-forms can be defined as $\left< \omega_1,\omega_2\right>=\int_M\omega_1\cdot\omega_2,$ where $\omega_1\cdot \omega_2$ is the sum of products of corresponding components of $\omega_1$ and $\omega_2$ expressed in an orthonormal frame. The codifferential operators $\delta_k: \Omega^k(M)\rightarrow \Omega^{k-1}(M)$ is defined as the adjoint operators of $d_{k-1}$ through $\left<\delta\omega_1,\omega_2\right>=\left<\omega_1,d\omega_2\right>$ for arbitrary $k$-form $\omega_1$ and $k-1$-form $\omega_2$. It is straightforward to verify that $\delta_{k-1} \circ\delta_k =0.$ The spaces of differential forms create a cochain complex or more specifically the de Rham complex,
\begin{center}
	\begin{tikzcd}
		\Omega^0(M) \arrow[r,yshift=1.7,"d^0"]
		& \arrow[l, yshift=-1.7, "\delta_1"] \Omega^1(M) \arrow[r,yshift=1.7,"d_1"]
		& \arrow[l, yshift=-1.7, "\delta_2"] \Omega^2(M) \arrow[r,yshift=1.7,"d_2"]
		& \arrow[l, yshift=-1.7, "\delta_3"] \Omega^3(M).
	\end{tikzcd}
\end{center}

The centerpiece of the de Rham-Hodge theory is the relation between the de Rham complex and the spaces of harmonic forms $\mathcal{H}_\Delta^k = \ker d_k \cap \ker \delta_k$ (which are naturally isomorphic to the corresponding cohomology groups $\ker d_k / \ima d_{k-1}$) as \[ \ker d_k = \mathcal{H}^k \oplus \ima d_{k-1}. \] Furthermore, \[\mathcal{H}^k(M)\cong \mathcal{H}_\Delta^k(M),\] where $\mathcal{H}_\Delta^k(M)=\{\omega|\Delta \omega = 0\}$ is the kernel of the Laplace-de Rham, or de Rham-Hodge Laplacian, operator $\Delta \equiv d\delta+\delta d=(d  + \delta )^2$. The equivalence with the above definition of $\mathcal{H}_\Delta^k(M)$ as the intersection of the kernels of $d$ and $\delta$ is due to $\left<\omega, \Delta \omega\right>=\left<d\omega,d\omega\right>+\left<\delta\omega,\delta \omega\right>$ for boundaryless $M$. Through Stokes' theorem, cohomology groups are isomorphic to homology groups, which allow for the categorization of surfaces. Specifically, Betti numbers, $\beta_k=\dim\mathcal{H}^k$, give the count for the number of $k$-dimensional holes and can be found by calculating the size of the null space of $\Delta$. Fundamentally, these harmonic forms are the 0-frequency spectral bases for differential forms, which allow us to study the structure of $M$. Instead of defining the Laplace operator with $\delta$, it can be rewritten using the Hodge star operator $\star$ which identifies a $k$-form with its complementary or dual $(n-k)$-form for an $n$-manifold by identifying the orthonormal basis of $k$-forms with the corresponding basis of $(n-k)$-forms, e.g., in 3D, the $2$-form basis $(dx^2\wedge dx^3, dx^3\wedge dx^1, dx^1\wedge dx^2)$ can be identified with $(\star dx^1, \star dx^2, \star dx^3),$ and the $3$-form basis $dx^1\wedge dx^2 \wedge dx^3$ can be identified with $\star 1$. With the Hodge star, the codifferential can be redefined as $\delta_k:=(-1)^k \star d\star$ (in 3D) so that the scalar Laplacian is $\Delta_0 = \star d\star d$ and the $k$-Laplacian is
\begin{equation}\label{eqn:kDelta}
    \Delta_k=(-1)^{k+1}\star d\star d+(-1)^k d\star d\star.
\end{equation}
Note that the $\Delta_1$ and $\Delta_2$ can be regarded as the vector Laplacian as they act on vector fields.

\subsection{Boundary conditions} \label{boundarycond}

Before further discussing the differential form version of the Laplacian on domains embedded in 3D Euclidean space, certain boundary conditions must be considered to ensure the operator is well-posed. Two common and useful choices are to have a differential form $\omega$ either normal or tangential to the boundary. In the case of normal boundary conditions, also known as (homogeneous) Dirichlet boundary conditions, $\omega$ is zero when applied to tangent vectors on the boundary. Using familiar calculus definitions for a scalar function $f$, the Dirichlet boundary condition is $f\vert_{\partial M}=0$. Similarly, for tangential boundary conditions, also called a Neumann boundary condition, $\star\omega$ is zero when applied to tangent vectors on the boundary. For instance, such a condition on a 1-form represented by a vector field $\bv$ (i.e., $\omega=v_1 dx^1+v_2 dx^2 + v_3 dx^3$)  can be expressed as
\begin{eqnarray*}
\star\omega(\bt,\cdot)&=&(v_1 dx^2\wedge dx^3+ v_2 dx^3\wedge dx^1 + v_3 dx^1\wedge dx^2)\left(t^1 \frac{\partial}{\partial x^1} + t^2 \frac{\partial}{\partial x^2}+t^3 \frac{\partial}{\partial x^3}\right) \\
&=&(v_2 t^3-v_3 t^2) dx^1+ (v_3 t^1-v_1 t^3) dx^2+(v_1 t^2-v_2 t^1) dx^3=0
\end{eqnarray*}
for any vector $\bt$ tangent to the boundary, or equivalently $\bn\cdot \bv\vert_{\partial M}=0,$ where $\bn$ is the boundary normal, i.e., the normal component of $\bv$ is $0.$

As a consequence, an update to the continuous setting is required when speaking in terms of vector fields. Enforcing that for tangential vector fields, i.e., normal 2-forms or tangential 1-forms, means satisfying one Dirichlet condition $\bv\!\cdot\!\bn\!=\!0$ and two Robin conditions $\bt_1\!\cdot\! \nabla_{\bn}\bv + \kappa_1 \bt_1 \!\cdot\! \bv\!=\!0, \bt_2\!\cdot\! \nabla_{\bn}\bv + \kappa_2 \bt_2 \!\cdot\! \bv\!=\!0$, here $\bt_1$ and $\bt_2$ are the two local tangent directions forming an orthonormal coordinate frame with surface normal $\bn$ and $\kappa_1$ and $\kappa_2$ are the sectional curvatures along the coordinate directions. These conditions effectively enforce no contribution from the curl along tangential directions. On the other hand for normal vector fields, i.e., tangential 2-forms or normal 1-forms, this means satisfying two Dirichlet conditions $\bv \!\cdot\! \bt_1\!=\!0$ and $\bv \!\cdot\! \bt_2\!=\!0$ and one Robin condition $\bn\!\cdot\!\nabla_{\bn} \bv+2H\bn\!\cdot\!\bv \!=\!0$, where $H$ is mean curvature, to enforce zero contribution from divergence.\footnote{These conditions are consistent with the differential form boundary conditions in Refs.\cite{zhao20193d,zhao2020rham,chenevol}. However, their vector field boundary conditions are correct only for flat boundary surfaces.}
Given these extra conditions, our harmonic space is now of finite dimension so that the kernel of the Laplacians is also finite dimensional and agrees with the corresponding absolute/relative homology dimensions.

In terms of differential forms, let $\Omega_t$ be the space of tangential forms with tangential differential so that $\omega_t \in \Omega_t$ if and only if $\star\omega_t\vert_{\partial M}=0$ and $\star d\omega_t\vert_{\partial M}=0$. One may interpret it that $\omega$ has ${n \choose k}$ DoFs per point, whereas $\star \omega\vert_{\partial M}$ has ${{n-1}\choose {n-k}}$ DoFs, so $\star d\omega\vert_{\partial M}$ provides the additional ${{n-1}\choose {n-(k+1)}}= {{n}\choose {k}}- {{n-1}\choose {n-k}}$ DoFs. In particular, for Neumann boundary condition of a 0-form $f$, $\star f\vert_{\partial M}=0$ is automatic since its an $n$-form evaluated in an $(n\!-\!1)$-D space, and $\star d f \vert_{\partial M}$ provides the familiar $\bn\cdot\nabla f\vert_{\partial M}=0.$ Similarly, let $\Omega_n$ be the space of normal forms with normal codifferential so that $\omega_n\in \Omega_n$ if and only if $\omega_n\vert_{\partial M}=0$ and $\delta\omega_n\vert_{\partial M}=0$.

Enforcing such modified boundary conditions means to restrict $\Delta$ to $\Omega_t$ or $\Omega_n$. The restricted harmonic forms are $\mathcal{H}^k_n=\mathcal{H}^k\cap\Omega_n$ and $\mathcal{H}^k_t=\mathcal{H}^k\cap\Omega_t$ which are associated with the kernels of $\Delta_{k,n}$ and $\Delta_{k,t}$ respectively. Due to the dimensionality of $k=\{0,1,2,3\}$ and since we are considering two types of boundary conditions, normal and tangential, there are eight Laplacian operators denoted $\Delta_{k,n}$ and $\Delta_{k,t}$.
However, using the duality between $k$ and $(n-k)$-forms, the space of normal $k$-forms can be identified with the space of tangential $(n-k)$-forms. This reduces the eight spectra to four distinct spectra, which may be further reduced to three using \emph{continuous} Hodge decomposition~\cite{zhao2020rham}.

\subsection{Discretization of the de Rham-Hodge theory}\label{sec:discrete}

To highlight the common algebraic structure underpinning both combinatorial Laplacian and de Rham-Hodge theory, we provide a combinatorial Laplacian theory perspective for the discretization  of the de Rham-Hodge theory.
 Kirchhoff introduced the burgeoning field of graph theory and the graph Laplacian matrix to circuit theory in 1847 \cite{kirchhoff}. Nearly 100 years later, Eckmann \cite{eckmann45} introduced   simplicial complexes to graphs, resulting in   higher-order graph Laplacians, also known as   combinatorial Laplacians. These tools and others in discrete calculus have been applied to topics firmly rooted in discrete data such as networks. For a further discussion of discrete calculus and discrete domains where   combinatorial Laplacians and their theory are applied, please see \cite{grady}.

When starting from this discrete, graph theoretic perspective, the input graph is a collection of objects and relationships between them denoted $G=(V,E)$ with vertex set $V=\{1,\dots , n\}$ and edge set $E\subseteq {V\choose 2}$.
Such relationships can be extended to higher orders using the language of graph cliques. A $k$-clique $K_k(G)$ is a set of $k$ vertices $\{v_1,v_2,\dots,v_k\}$ such that $e_{ij}=\{v_i,v_j\}\in E$ for every $\{v_i,v_j\}\in K_k(G)$. Similarly, a $k$-simplex can be discretely defined as an ordered set of $k+1$ vertices, $\sigma_k=[v_0,v_1,\dots,v_k]$. For example, triangles can be denoted by their vertices as 2-simplices. An abstract simplicial complex $K$ is thus an extension of graph to include simplices with dimensions higher than 1, provided that any subset of $k$-simplex is also included in the complex, i.e., $\sigma_j\in K$ and $\sigma_k\subseteq \sigma_j$ then $\sigma_k \in K$. 
Using cliques, an abstract simplicial complex can be formed by a collection of such clique sets since the above condition is satisfied. This simplicial complex is called the clique complex, and any other abstract simplicial complex with edges (1-simplices) given by the graph $G$ can be seen as a subcomplex of it. For instance, one commonly used simplicial complex called Vietoris-Rips complex is a special case of clique complex. One can further specify particular functions defined on vertices, edges and higher order cliques as the discrete analogs of differential forms on manifolds, also called cochains.


As an important special case for concrete computations, these definitions can be given different specifications for working with graph-like objects embedded in $\mathbb{R}^n$ where a $k$-simplex $\sigma_k$ is defined instead as the convex hull of $k+1$ affinely independent points in $\mathbb{R}^n$.
For our purposes, $k\in\{0, 1, 2, 3\}$, so $k$-simplices correspond to their $k$-dimensional elements: vertices, edges, faces, and cells. It is useful to use the terms faces and cofaces when referring to relations between simplices. A $k$-simplex can be referred to as a $k$-face of a simplex $\tau$, if $\sigma$ is on the boundary of $\tau$, i.e., it is a simplex formed from a subset of the vertices of $\tau$. In this case, $\tau$ is a coface of $\sigma.$
A simplicial complex can be seen as a geometric realization of an abstract simplicial complex with vertices identified as points in $\mathbb{R}^n$ such that the non-empty intersection of any two simplices $\sigma_j, \sigma_k \in K$ is a face of $\sigma_j$ and $\sigma_k$. Furthermore, discrete differential forms in this context can be defined as  functions on these simplices correspond to integrated geometric quantities, such as electric field, magnetic flux, and mass density. The concepts can be generalized to manifolds, which, roughly speaking, are spaces that are locally $\mathbb{R}^n$.

Such assumptions are useful when discretizing some presumed continuous shape for application purposes. As shown in the Venn diagram in Figure \ref{fig:venn}, these assumptions lead to a discretized Hodge Laplacian, related to but not to be confused with the combinatorial Laplacian (we use the term in the sense that they are Laplacians defined purely on simplicial complex,  without weights or functions assigned to simplices, as in \cite{friedman}). While both combinatorial and discretized Hodge Laplacian are categorized as discrete Laplacians, the discretized Hodge Laplacian is derived from the usual Hodge Laplacian from the continuous theory.
The discrete treatment of de Rham-Hodge theory, which preserves many properties of its continuous counterpart, can be described using the language of differential forms as acted on by boundary and Hodge star operators by  DEC \cite{dec06}.
%

To discretize de Rham-Hodge theory, first the input domain must be specified. It is common to collect sample points of the input shape, where each point is presented by its 3D coordinates. The relationship among points may be given through mesh connectivity (cell complex structure) or prescribed by any number of criteria. A Delaunay triangulation~\cite{delaunay1934sphere} for a given set of vertices is defined by its dual structure, the corresponding Voronoi diagram~\cite{voronoi1908nouvelles}. A Voronoi dual cell of a vertex $v$ is the set of points $P$ around $v$ such that the distance between any $p\in P$ is less than the distance between $p$ and any other vertex,
\begin{equation*}
    V_i = \{v\in \mathbb{R}^q \ | \|v-p_i\| \le \|v-p_j\|, \text{ for all } p_j \in P \}.
\end{equation*}
The Voronoi diagram is the set of Voronoi cells, which is defined as
\begin{equation*}
    \text{Vor}P = \{V_i \ | \ \text{ for all } i \in \{1,2,\cdots,|P|\}\},
\end{equation*}
where $|P|$ is the number of points in set $P$.
The primal 2D (3D, respectively) Delaunay triangulation can also be defined by the empty circle (sphere) property, that the circumscribing circle (sphere) of a face (cell) does not contain any vertex in its interior (cell). This triangulation is unique as long as any four vertices are not co-circular (or any five vertices are not co-spherical).

Given a discrete dataset, different simplicial complexes can be defined with simple criteria such as the alpha complex, \v{C}ech complex, and Vietoris-Rips complex. Given a set of vertices $V =\{v_0, v_1, \cdots, v_{N_0-1}\}$, embedded in $\mathbb{R}^3$, consider a nested family of simplicial complexes. This family, called a filtration, is created for a positive real, integer number $\alpha$, the filtration parameter. Specifically, using the alpha complex filtration as an example, the filtration of subcomplexes $(K_\alpha)_{\alpha=0}^m$ is
\begin{equation*}
    \emptyset = K_0 \subseteq K_1 \subseteq K_2 \subseteq \cdots \subseteq K_m = K.
\end{equation*}
which gives the final simplicial complex $K$, the Delaunay triangulation. Each subcomplex can be considered as an $\alpha$-complex, the collection of Delaunay simplices whose smallest empty circumsphere has a radius $\le\alpha$. Such simplicial complexes are computed in the field of TDA in the context of persistence of features over time where the persistent combinatorial Laplacian is computed for varying filtration parameter values  \cite{wanggraph}.

The boundary operator $\partial$ of a $k$-simplex is a $(k-1)$-chain (i.e., a formal linear combination of  $(k-1)$-simplices)
\begin{equation}\label{eqn:boundaryop}
    \partial \sigma=\sum^k_{i=0}(-1)^i[v_0,v_1,\dots,\hat{v_i}\dots,v_k],
\end{equation}
where $\hat{v_i}$ is the vertex omitted. It can be linearly extended to a linear operator on $k$-chains. In the matrix form, the $k$-boundary operator $B_k$ has the number of $(k-1)$-simplices rows and the number of $k$-simplices columns. In the matrix, an $(i,j)$ entry is $+1$ or $-1$ if the $i$-th $(k-1)$-simplex is on the boundary of the $j$-th $k$-simplex and 0 if it is not on the boundary. The sign is determined by the orientation of the $(k-1)$-simplex relative to the orientation of the $k$-simplex. The orientations can be arbitrarily chosen for the simplices forming a basis for the space of chains. Next, since the boundary operator $\partial$ is the adjoint of the exterior derivative $d$, it is also called the coboundary operator. One can take the transpose of the proposed $B_k$ to get $D_{k-1}$, the discrete representation of $d$ when $k$-forms $\omega$ are discretized to $k$-cochains $w$, the dual of $k$-chains, by $w_\sigma=\int_\sigma \omega$. This equivalence to the boundary operator can also be seen through Stokes' theorem,
\begin{equation*}
    \int_{\partial\sigma}\omega=\int_\sigma d\omega,
\end{equation*}
which means that calculating $d\omega$ on a simplex $\sigma$ can be replaced with a calculation of $\omega$ on the boundary of $\sigma$. Now $D_k = B_{k+1}^T$ is the signed incidence matrix between $k$ and $(k+1)$-simplices.

To discretize the Hodge star $\star$, the concept of duality between Delaunay triangulation and Voronoi diagram is further explored. Since the Hodge star maps $k$-forms to $(n-k)$-forms, we need a way to transform primal $k$-simplices to dual $(n-k)$-cells and vice versa. One simple way to calculate the Hodge star, is to discretize $\star \omega$ by its integral on Voronoi dual cells when $\omega$ is represented by the integral on the Delaunay simplices. In this case, we call $\omega$ as a primal form and $\star \omega$ a dual form. The Hodge star can then be discretized as a diagonal matrix $S_k$ as the ratio between the   length, area, or volume of a dual $(n-k)$-Voronoi-cell and its primal $k$-simplex. For this reason, the Hodge star is essentially a scaling factor and can be highly dependent on the tessellation. Other more accurate Hodge stars can be computed, such as the (non-diagonal) Galerkin Hodge star \cite{Bossavit2000ComputationalEA} through higher-order basis functions. Note that the Hodge stars are not necessarily unique even within the same computational domain, since they are determined by measurements of cells that are influenced by the material properties relevant to the actual applications. Moreover, both $D$ and $S$ operators can be similarly defined on cell complexes, a generalization of simplicial complexes. Cell complexes can approximate manifold domains by tessellating them into topological ball-like cells, whose boundaries are unions of lower dimensional cells.

Next, the discretized Hodge Laplacian follows from \ref{eqn:kDelta},
\begin{equation}\label{eqn:hodgelaplace}
    L^H_k=D_k^TS_{k+1}D_k+S_kD_{k-1}S_{k-1}^{-1}D_{k-1}^TS_k,
\end{equation}
which is a symmetric matrix where $S_k^{-1}{L}^H_k$ is the discrete counterpart of $\Delta_k$.
The signs in the above equation are consistent with the proper use of $D_k^T$ to discrete $(-1)^k d_{n-k-1}.$
As discussed, to discover geometrical information, it helps to investigate not only the null space of the Laplacian but the entire spectra, the eigenvalues of the Laplacian. While the size of the null space of the $k$th-Laplacian, the $k$th-Betti number, reveals the number of $k$-dimensional holes, we would also like to investigate \emph{non-zero eigenvalues} which offer rich qualitative information of the shape, especially in the context of machine learning. In fact, the Fiedler value, i.e., the first non-zero eigenvalue, can describe connectivity. Furthermore, the multiplicity of eigenvalues may reveal certain symmetries of a shape.

\subsection{Combinatorial Laplacian}\label{sec:cliquegraph}

 In an apparent analogy to the Hodge Laplacian,  the combinatorial Laplacian has a similar algebraic structure and is defined as
\[ L^G_k=D_k^TD_k+D_{k-1}D_{k-1}^T.\]
Then $L^G_0$ is called the graph Laplacian as used in spectral graph theory~\cite{lim15}. Here Betti numbers correspond to the dimension of nontrivial $k$-dimensional cycles (i.e., boundaryless $k$-chains that are not boundaries of $(k+1)$-chains). Simply stated, the combinatorial Laplacian disregards the geometric measurements  as encoded by the Hodge star, e.g., edge length, face area, or cell volume. As discussed in  \cite{lim15}, the combinatorial Laplacian only requires the adjacency information contained in boundary operators. Furthermore, the discrete operators proposed are sparse and positive semi-definite, which allows for the use of efficient solvers. However, when the connectivity information given is restricted to a graph, while the $0$-Laplacian can incorporate edge weights to have a weighted combinatorial Laplacian, there is no generic way to assign proper weights to higher-order simplices other than the actual Hodge star calculation, which are almost never a scaled identity matrix.

For a clique complex, the resulting $k$-Laplacian is almost never that of the (graph) Laplacian of a simplicial mesh. \v{C}ech complexes also do not share the same $D_k$ with a simplicial mesh. While the alpha complex does provide incidence matrices similar to a manifold mesh aside from some degeneracy, the resulting Laplacians without $S_k$ still produce spectra that are far from those of  Hodge Laplacians.  Recall from graph theory that the degree of a vertex $v$ is the number of adjacent vertices, denoted $\deg(v)$. With the clique complex, the adjacency information between vertices induces incidence, i.e., the face/coface relationship between $k$ and $(k\!+\!1)$-simplices. For example, given the end point vertex $v$ of an edge $e$, $v$ is incident to $e$.

Wang et al.~\cite{wanggraph} discussed the general case when the simplicial complex is not necessarily induced from the clique in the graph. It is useful to introduce the notion of adjacency, incidence, and degree for all simplices in the complex. Two simplices $\sigma^i$ and $\sigma^j$ are called upper adjacent, denoted $\sigma_k^i\overset{U}{\sim}\sigma_k^i$, if they are incident to the same $(k+1)$-simplex. The upper degree $\deg_U(\sigma_k)$, is the number of $(k+1)$-simplices that are faces of $\sigma_k$. Similarly, two $k$-simplices are called lower adjacent if they are adjacent, or share the same $(k-1)$-simplex. The lower degree $\deg_L(\sigma_k)$, is the number of $(k-1)$-simplices that are cofaces of $\sigma_k$. Then the degree of a $k$-simplex, $k>0$ is,
\[\deg(\sigma_k)=\deg_U(\sigma_k)+\deg_L(\sigma_k).\]
To explicitly define the combinatorial Laplacian, an orientation for the simplicial complex must be specified. While the orientation itself is arbitrary, each simplex, except for vertices, is given a direction. For simplex $\sigma_k$, its ordering is defined by an ordering on its vertex set, $\sigma_k=[v_0,\dots, v_k]$. The combinatorial Laplacian matrix is an $N_k\!\times N_k$ matrix, where $N_k$ is the number of $k$-dimensional simplices. Its entries are explicitly defined as follows\footnote{We follow the notations in  \cite{wanggraph}, except for their typo in the diagonal entry, where an addition should be replaced with an equality, that is, when $i=j$,  $(L^G_k)_{ij}=\deg(\sigma^i_k) = k+1$.},

\begin{equation}\nonumber
\left(L^G_0\right)_{ij}=
\begin{cases}
	\deg(\sigma_0^i) & \text{ if } i=j. \\
	-1 & \text{ if } \sigma_0^i\overset{U}{\sim}\sigma_0^j.\\
	0 & \text{otherwise},
\end{cases}
\end{equation}
and for $k>0$
\begin{equation} \nonumber
\left(L^G_k\right)_{ij}=
\begin{cases}
	\deg(\sigma_k^i) & \text{ if } i=j. \\
	1 & \text{ if } i\neq j,\sigma_k^i\overset{U}{\nsim}\sigma_k^j\text{ and }\sigma_k^i\overset{L}{\sim}\sigma_k^j \text{ with the same orientation.}\\
	-1 & \text{ if } i\neq j,\sigma_k^i\overset{U}{\nsim}\sigma_k^j\text{ and }\sigma_k^i\overset{L}{\sim}\sigma_k^j \text{ with different orientation.}\\
	0 & \text{ if } i\neq j \text{ and either }\sigma_k^i\overset{U}{\sim}\sigma_k^j\text{ or }\sigma_k^i\overset{L}{\nsim}\sigma_k^j.\\
\end{cases}
\end{equation}

The combinatorial Laplacian  can certainly give rise to its own decomposition theory, similar to the inner products assigned to the cochain spaces (including replacing the Hodge star by the identity matrix) lead to a corresponding Hodge theory~\cite{eckmann45}. Discrete Hodge operators can draw analogy from their continuous counterpart. For instance, in electromagnetics, Kirchhoff's laws on electrical circuits can be regarded as the examples of discrete curl applied to electric fields and discrete divergence applied to currents, which are indeed direct consequences of the continuous theory. Note that only carefully chosen discrete Hodge stars lead to actual approximation of the Hodge theory of \emph{differential} forms~\cite{dodziuk}, e.g., for computing eigenvalues in the spectra. However,  combinatorial Laplacians do not correspond to a continuous theory in general.

\section{Boundary-Induced Graph (BIG) Laplacian}\label{sec:big}

\begin{figure}[H]
\label{fig:decomposition}
\centering
\includegraphics[width=.9\textwidth]{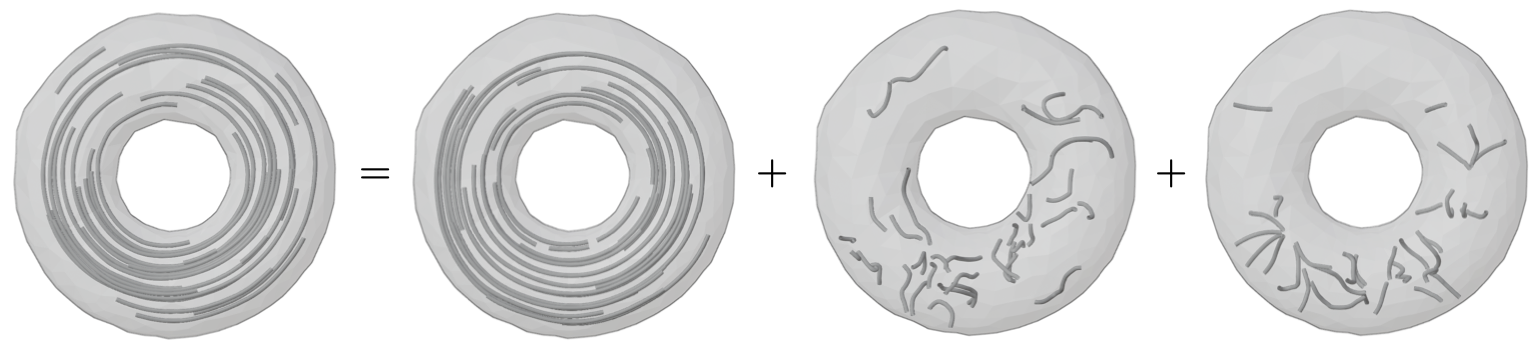}
\caption{Hodge decomposition of a vector field with zero combinatorial Laplacian. Given an irregularly sampled volumetric mesh of torus (with one half more densely sampled), we used the combinatorial Laplacian to generate a combinatorially harmonic field (left). We then computed a decomposition using the discretized Hodge Laplacians to compute the harmonic component (middle left), a divergence-free component (middle right), and a curl-free component (right). Due to the simplified Hodge star in the combinatorial Laplacian, the \emph{combinatorially} harmonic field (left) is indeed neither curl-free nore divergence-free as seen in the right two components.}
\end{figure}

While the combinatorial Laplacian based on cliques is useful for calculating Betti numbers, it generally differs drastically from the Hodge Laplacian  and the discretized Hodge Laplacian. This is due to the fact that the combinatorial Laplacian is fully determined by graph connectivity and lacks geometry awareness. Only if data points are uniformly distributed with regular graph edges can the spectra be comparable given perhaps some scaling factors. However, such a distribution is only possible on a 2D plane given explicit graph and clique complex representations. Furthermore, physical applications require some notion of behavior on the boundary, not generally available in the context of the combinatorial Laplacian. Even for domains without boundary, the combinatorial Laplacian can lead to results that are unrecognizable in contrast to those generated with geometric measurements and manifold-based connectivity. Hodge Laplacians based on the connectivity from a quadrilateral mesh can be used to solve for the decomposition of an input vector field into the sum of a gradient field, a rotational field, and a harmonic field, where the Betti numbers are $(\beta_0,\beta_1,\beta_2)=(1,2,1)$. In contrast, for the clique-based combinatorial Laplacian built on the 1D skeleton graph of the quad mesh, there are no 2D cells (no 3-cliques) for the decomposition, and even the Betti numbers are completely different $(\beta_0,\beta_1,\beta_2)=(1,1999,0)$ given that the $2000$ original quads are missing. Of course, combinatorial Laplacians based on cell complexes preserving the topology would lead to the right Betti numbers. However, even in such cases, comparing the decomposition of the discretized Hodge Laplacian to the combinatorial Laplacian on a Delaunay tetrahedralization in Fig.~\ref{fig:decomposition}, only the results based on discretized Hodge Laplacian approximate the continuous, physical meaning of decomposition into curl and divergence-free vector fields. Nevertheless, it is possible to have both the simplicity of combinatorial Laplacians and the approximation to the continuous Laplacian through minimal changes to the combinatorial Laplacians for entries associated with boundary cells, as we show though our novel 3D boundary-induced graph  Laplacian defined implicitly on a regular Cartesian grid.


\subsection{Eulerian representation}

From the context of applications such as fluid simulation, the distinction between explicit and implicit methods can be described by Lagrangian and Eulerian discretizations. The Lagrangian method uses cell complexes built on vertices that are irregularly sampled such as the fluid particles/parcels at a given time, while the Eulerian representation uses vertices fixed in the domain. As commonly used in level set-based methods, we employ an Eulerian method using a Signed Distance Function (SDF) to the boundary surface of a volume $M$ on a regular Cartesian grid as opposed. This allows the underlying connectivity to be determined without an explicit incidence matrix representation, and the Hodge stars to be almost rescaled identity matrices, while ensuring the spectra of BIG Laplacians to resemble those of the continuous Hodge Laplacians on the volumed bounded by a (piecewise-)smooth surface.

In addition to facilitating the proper comparison of the different Laplacians, a regular grid also helps to simplify the data structures used. Since grid elements have a fixed length, area, and volume, the effect of the Hodge star, or lack thereof in the case of the BIG Laplacian, is minimized. The grid spacing or grid length is given as input and though not proven in this work, a reduction in grid length is empirically shown to increase the accuracy of our method as demonstrated in Section \ref{sec:sim}. Similar to the simplicial complex differential operator, our $D_k$ includes grid elements forming a cell complex instead of a simplicial complex.

\subsection{Boundary condition correction} \label{NTC}

As already included in the study of the usual discrete Hodge Laplacian, boundary conditions are now introduced for the BIG Laplacian. By simply checking the sign of SDF values of grid points, we can tell which cells are on the boundary of $M$ and thus enforce particular boundary conditions through updates to $D_k$. As discussed, two types of boundary conditions are considered, normal (Dirichlet) and tangential (Neumann).  In the discrete Lagrangian setting (simplicial complexes, i.e., triangle/tetrahedron meshes in 2D/3D), these conditions correspond to exclusion, in the normal case, or inclusion, in the tangential case, of boundary elements. The exclusion is done directly in $D_k$ operators as rows and columns on the boundary are kept or removed. Accordingly, we denote the differential operator for normal forms as $D_{k,n}$ and for tangential forms as $D_{k,t}$.

In our Eulerian setting, the cells are not aligned with the boundary as we do not have an explicit boundary surface to enable the exclusion of the surface to the volume. Fortunately, on regular grids, both boundary conditions can be implemented based on the cells of either the primal or dual grid. For normal boundary conditions, we include all cells that are either inside or intersect with the boundary., i.e., at least one of its vertices is inside. For tangential boundary conditions, we include all cells whose corresponding dual cells are inside or partially inside. In fact, since the dual grid structure is also a Cartesian grid staggered with the primal grid (by a displacement of $\frac12(l_g,l_g,l_g)$) we only need to implement one of the boundary conditions $L_{k,n}$ as the other $L_{n-k,t}$ would be strictly equivalent to a slightly shifted SDF input (or equivalently, a slightly shifted grid).

Then, to facilitate the matrix conversions, a projection or selection matrix $P_k$ may be used. $P_k$ removes any columns or rows corresponding to $k$ and $k+1$ grid elements to be removed as either exterior to the model or as prescribed by the boundary condition. An example is provided in Sec.~\ref{sec:example2DNormalBC}.
Using the projection matrices, the new exterior matrices can be expressed as
\begin{equation*}
    D_{k,n} = P_{k+1}D_kP_{k}^T,
\end{equation*}
and the normal BIG Laplacians are
\begin{eqnarray}\label{eq:normalGraphLaplacians}
    L^B_{0,n} &=& D_{0,n}^TD_{0,n} \\
    L^B_{1,n} &=& D_{1,n}^TD_{1,n}+D_{0,n}D_{0,n}^T \\
    L^B_{2,n} &=& D_{2,n}^TD_{2,n}+D_{1,n}D_{1,n}^T \\
    L^B_{3,n} &=& D_{2,n}D_{2,n}^T.
\end{eqnarray}

To solve for the spectra, we use a scaling diagonal matrix, so that the results of the BIG Laplacian are comparable to the Hodge Laplacian as calculated in the next section. In particular, let $S^B_k$ have $1/l_g^2$ on the diagonal with correct size for the number of $k$ grid elements.
Then with projection matrices, $P_k(S^B_k)^{-1}P_k^TL_k$ corresponds to $\Delta_k$.

\subsection{Hodge Laplacian on regular grids}

Now the BIG Laplacian can be directly compared to the usual Hodge Laplacian, though we include some updates to the Hodge Laplacian on the regular grid for completeness. For simplicity and numerical reasons, the $k$-Hodge star is usually implemented as a diagonal matrix $S_k$ with entries as the ratio between the $(n-k)$-volume of the Voronoi (or other dual structure) $(n-k)$-cell and the $k$-volumes of the primal $k$-simplex and is highly dependent on the tessellation. More accurate Hodge stars can be discussed but are left for future work on numerical accuracy. Yet, given our regular grid representation, Hodge stars are simple and uniform to calculate for grid elements located completely in the interior of the model. With grid edge length $l_g$, values of $S_0$ are $l_g^3$, $S_1$ are $l_g^2/l_g$, $S_2$ are $l_g/l_g^2$, and $S_3$ are $1/l_g^3$. However, we modify the Hodge star slightly for edges, faces, and cells that cross the boundary surface following \cite{liu2015fluid} with extensions to all orders for both boundary conditions, since they only implemented normal $1$-forms representing vorticity fields of tangential velocity fields.

As with $D_k$, boundary conditions should also be considered for $S_k$. In particular, For $S_{k,n}$, the primal cell volumes are modified with dual cell volume unaltered, whereas for $S_{k,t}$, dual cell volumes are modified with the primal cell volume unaltered. If the grid is shifted by half a grid edge length in all three directions, the new $S_{k,t}$ has nonzero entries that are the inverse of the corresponding entries in the original $S_{3-k,n}.$ For instance, in the case of 2-forms $\omega$ with normal boundary conditions, in order to define $S_{1,n}$, $S_{2,n}$, and $S_{3,n}$, we require that $\omega$ is zero when applied to tangent vectors on the boundary. This is equivalent to saying that surface patches in the Lagrangian sense offer no contribution on the boundary or that the flux through the boundary is zero. As before with $D_{k,n}$ we use the projection matrices $P_k$ as follows,
\begin{equation*}
    S_{k,n} = P_{k}S_kP_k^T.
\end{equation*}

For example, we alter $S_2$ so that grid faces which cross the boundary only contribute as much as the portion of the face which resides inside of the model. For the normal boundary condition, the dual cell sizes are not altered \cite{liu2015fluid}. Thus, the Hodge star $S_2$ acts as a scaling factor, a ratio of original dual edge length to the portion of primal face area inside of the model $\mathcal{A}_i$. Given that the grid spacing is $l_g$,
\[ S_{2,n}(i, i) = \frac{l_g}{\mathcal{A}_i} \]
for face $f_i$ on the interior or crossing the boundary of $M$. Note that for faces $f_{j}$ completely inside of the model, $S_{2,n}(j, j)=1/l_g$. The situation is similar for $S_{1,n}$ and $S_{3,n}$ in that the primal length or volume inside the model is considered. Approximations of the inside portions are computed using the SDF values, the convex hull of inside simplices, and the marching cubes algorithms. Here, $S_0$ is unchanged since a vertex is dimensionless, so it is either inside or outside. For the numerical implementation, we avoid division by near zero primal cell volumes. Any primal $k$-cell with a volume below a preset positive threshold ($\epsilon\ll l_g$) will have its volume rounded up to $\epsilon^k$. This prevents the large scaling factors that can potentially distort the results, while maintaining the correct Laplacian kernel dimensions.

The normal Hodge Laplacian uses the $D_k$ operators as in the BIG Laplacian but with $S_k$ operators for scale and are
\begin{eqnarray}\label{eq:normalLaplacians}
    L^H_{0,n} &=& D_{0,n}^TS_{1,n}D_{0,n} \\
    L^H_{1,n} &=& D_{1,n}^TS_{2,n}D_{1,n}+S_{1,n}D_{0,n}S_{0,n}^{-1}D_{0,n}^TS_{1,n} \\
    L^H_{2,n} &=& D_{2,n}^TS_{3,n}D_{2,n}+S_{2,n}D_{1,n}S_{1,n}^{-1}D_{1,n}^TS_{2,n} \\
    L^H_{3,n} &=& S_{3,n}D_{2,n}S_{2,n}^{-1}D_{2,n}^TS_{3,n}.
\end{eqnarray}
All matrices are symmetric and positive semi-definite where $(S_{k,n})^{-1}L^H_{k,n}$ corresponds to $\Delta_k$. We solve for the eigenvalues and eigenbasis functions through a generalized eigenvalue problem,
\[ L^H_{k,n}\omega_k = \lambda_kS_{k,n}\omega_k.\] See algorithm \ref{alg:buildops} for more details.

\begin{algorithm}[ht!]
\caption{Assemble coboundary, Hodge star, and projection matrices.}
\label{alg:buildops}
\begin{algorithmic}[1]
\REQUIRE{dimension $k\in\{1,2,3\}$ \\
\quad\quad primal $k$-areas $\mathcal{A}_{k}^{\text{primal}}(i)$ for $i\in\{0,1,\dots,\vert\Sigma_k\vert-1\}$ \\
\quad\quad error tolerance $\epsilon\ll l_g$ \COMMENT{prevents division by $0$}}
\ENSURE{$D_{k-1}$, $S_k$, and $P_k$}
\STATE{$n\leftarrow0$}
\STATE{$P_k, S_k\leftarrow$ empty sparse $\vert\Sigma_{k}\vert\times\vert\Sigma_{k}\vert$ matrix}
\COMMENT{$\Sigma_k$: set of oriented $k$-cells}
\STATE{$D_{k-1}\leftarrow$ empty sparse $\vert\Sigma_{k}\vert\times\vert\Sigma_{k-1}\vert$ matrix}
\FOR{$\sigma_i \in \Sigma_k$}
\IF{ any incident gridpoint of $\mathcal{A}_{k,i}^{\text{primal}}$ is inside}
\STATE{$S_k(i,i) \leftarrow \mathcal{A}_{3-k}^{\text{dual}}/\max(\mathcal{A}_{k}^{\text{primal}}(i),\epsilon^k)$}
\COMMENT{dual $(3\!-\!k)$-area $\mathcal{A}_{3-k}^{\text{dual}}=l_g^{3-k}$}
\STATE{$P_k(n,i)=1$}
\STATE{$n\leftarrow n+1$}
\ENDIF
\FOR{$\xi_j\in\Sigma_{k-1}$ a coface of $\sigma_i$}
\STATE{$D_k(i,j) = 1 \text{ or } -1$ \COMMENT{depending on orientation}}
\ENDFOR
\ENDFOR
\end{algorithmic}
\end{algorithm}

\subsection{The reduced spectrum $N, T$, and $C$}

Since the space of normal $k$-forms can be identified with the space of tangential $(n-k)$-forms the four independent operators in \eqref{eq:normalLaplacians} are sufficient to represent all eight Laplacians. Furthermore, as the Laplacians are constructed based on the three exterior derivatives, the spectra of these can be decomposed into three distinct parts~\cite{zhao20193d}:
\begin{enumerate}
    \item $T$, singular values of the gradient of tangential scalar fields (or equivalently, divergence of tangential gradient fields),
    \item $N$, singular values of the gradient of normal scalar fields (or equivalently, divergence of normal gradient fields),
    \item $C$, singular values of the curl of tangential curl fields (or equivalently, curl of normal curl fields),
\end{enumerate}
The spectrum of each Laplacian can be represented as a combination of one or two of the three singular spectra, potentially adding a zero with a multiplicity based on the corresponding Betti number. This can be illustrated by the fact that variations of $D_kD_k^T$ and $D_k^TD_k$ both show up in $L_k$ and $L_{k+1}$ and each has the same set of nonzero eigenvalues. For example, the spectrum of $L_{0,n}$ appears as part of the spectrum of $L_{1,n}$ from the gradient fields.

\section{Differences among topological Laplacians}\label{sec:diff}

To illustrate the differences among the usual combinatorial Laplacian, BIG Laplacian, and Hodge Laplacian, we introduce some simple examples discussing different boundary conditions and spatial discretizations. We start by calculating the Eulerian BIG and Hodge 0-Laplacians and their spectra for 2D shapes over coarse grids, as well as an example for a triangulated (Lagrangian) shape. Next, we show the drastic differences between the Lagrangian clique-based combinatorial and Hodge Laplacians in the case of an irregularly sampled ball. We conclude with a figure describing the importance of different boundary conditions and the effect it has on the spectra of topological Laplacians.

\subsection{2D calculation with normal boundary condition}\label{sec:example2DNormalBC}

\begin{figure}[H]\label{fig:example}
\centering
\includegraphics[width=.5\textwidth]{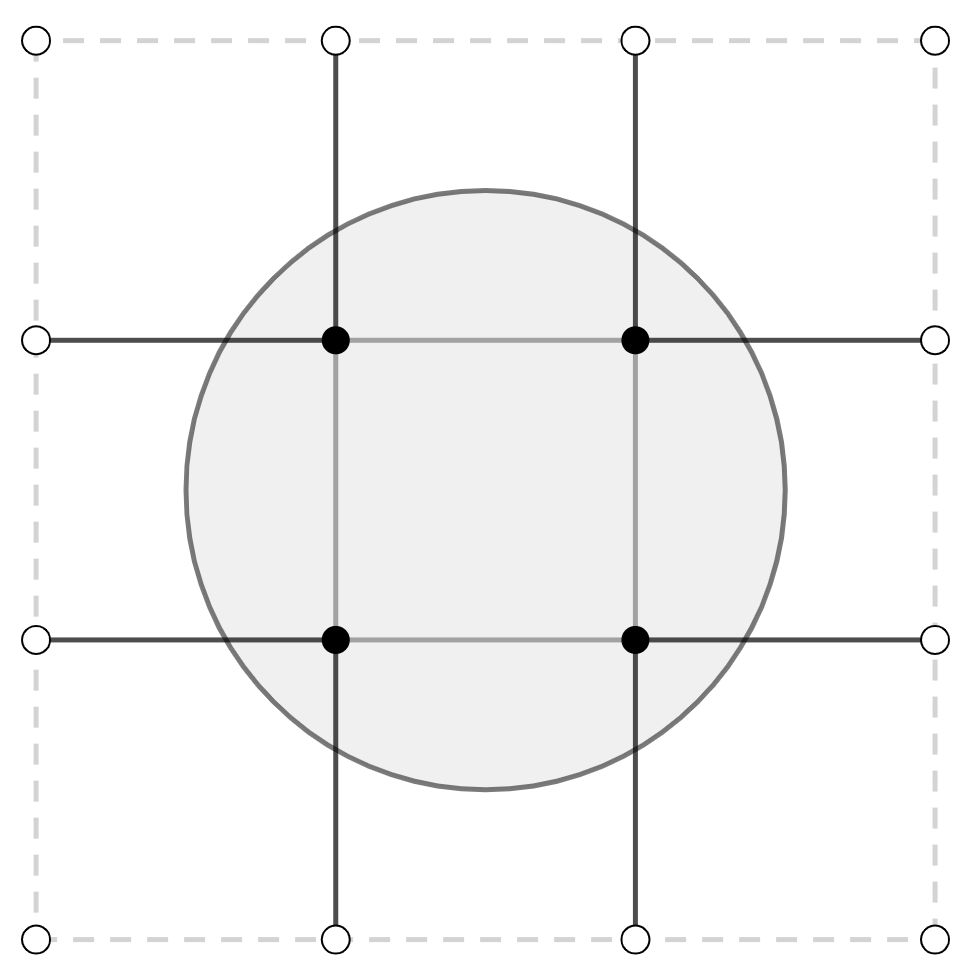}
\caption{A unit disk embedded in a $3\!\times\!3$ 2D grid with unit grid length.}
\end{figure}

In Figure \ref{fig:example}, we first present the necessary matrices to compute $L_{0,n}$ for two simple examples, a disk and a square on a $3\!\times\!3$ grid with grid length $l_g=1.$ The disk example illustrates the effect of the Hodge star on the Hodge Laplacian. The 8 edges which cross the boundary of the model (in black in the figure) contribute only their primal edge length, a fraction of the grid edge length to the calculation. In particular, $S_{1,n}$ is the ratio of the dual edge length to the (partial) primal edge length.

Given the unit disk with an embedding as in Figure \ref{fig:example}, to compute the 0-Laplacian with normal boundary conditions we  first remove the outside edges and vertices, the dashed gray edges and empty circle vertices in Figure \ref{fig:example}. For example, the removal of vertices can be regarded as the projection $P_0,$ whose entries are all $0$ except $P_{0,[1,6]}=P_{0,[2,7]}=P_{0,[3,10]}=P_{0,[4,11]}=1.$

The remaining edges are then split into two sets, the solid gray edges inside $E_i$ and the black edges which cross the boundary $E_b$. Then the coboundary operator is
\begin{equation}\label{eq:rectD0}
D_{0,n} =
\begin{bmatrix}
-1 &  1 &  0 &  0  \\
 0 & -1 &  1 &  0  \\
 0 &  0 & -1 &  1  \\
 1 &  0 &  0 & -1  \\
-1 &  0 &  0 &  0 \\
 0 &  1 &  0 &  0 \\
 0 & -1 &  0 &  0 \\
 0 &  0 &  1 &  0 \\
 0 &  0 & -1 &  0 \\
 0 &  0 &  0 &  1 \\
 0 &  0 &  0 & -1 \\
 1 &  0 &  0 &  0
\end{bmatrix}
\end{equation}
noting that the first four edge rows are for those in $E_i$ and have both $-1$ and $1$ entries for their vertex boundaries while the rest of the edges in $E_b$ only have one vertex included. Denote the primal edge length of edge $e_i$ as $l_{p}(i)$ then compute the diagonal Hodge star,
\[S_{1,n}(i,j) = 1/l_p(i)\;\delta_{ij},\]
where $\delta_{ij}$ is the Kronecker delta. Note that edges which cross the boundary have primal inside edge length $l_p =0.37$.

Next we compute the Laplacians,
\[L_{0,n}^H = D_{0,n}'^TS_{1,n} D_{0,n} =
\begin{bmatrix}
    7.41  & -1   &  0  & -1 \\
   -1  & 7.41  & -1 &  0 \\
    0  & -1   & 7.41  & -1 \\
   -1 &  0  & -1  &  7.41
\end{bmatrix},\]
\[L_{0,n}^B = D_{0,n}^T D_{0,n} =
\begin{bmatrix}
    4 & -1 &  0 & -1 \\
   -1 &  4 & -1 &  0 \\
    0 & -1 &  4 & -1 \\
   -1 &  0 & -1 &  4
\end{bmatrix}.\]
The eigenvalues of these matrices are then
\[\lambda_{0,n}^H=\{5.41, 7.41, 7.41, 9.41\}\text{ and }\lambda_{0,n}^B=\{2, 4, 4, 6\}.\]
While the exact eigenvalues for the Hodge Laplacian, as computed from the roots of the Bessel functions squared are
\[\lambda_\text{exact}=\{ 5.7832, 14.682, 14.682, 26.3746\}.\]
Thus, the Hodge Laplacian's principal eigenvalue is much closer than the BIG Laplacian to the expected eigenvalue even for a very coarse sampling of points thanks to the Hodge star. In fact, the clique-based combinatorial Laplacian produces an even more erroneous principal eigenvalue $0$, since it can only handle the tangential boundary condition as discussed in Section~\ref{sec:tangentialRect}.

For a $3\!\times\!3$ square, $D_0$ remains the same and $S_1$ is the $12\!\times\!12$ identity matrix since all 12 edges are entirely inside of the model. For this reason, we have
\[L_{0,n}^H = L_{0,n}^B =
\begin{bmatrix}
    4 & -1 &  0 & -1 \\
   -1 &  4 & -1 &  0 \\
    0 & -1 &  4 & -1 \\
   -1 &  0 & -1 &  4
\end{bmatrix},\]
\[\lambda_{0,n}^H= \lambda_{0,n}^B = \{2, 4, 4, 6\}, \text{ and for the Hodge Laplacian,}\]
\[\lambda_\text{exact}=\{ 2.1932, 5.4831, 5.4831, 8.773\}.\]

In this case, since the square's boundary is exactly the boundary of the grid domain, the Hodge star is the identity matrix. Thus, the Hodge and BIG Laplacians are identical. However, the eigenvalues in this case still perform favorably even for such a coarse tessellation. The accuracy of eigenvalues improves considerably with a finer resolution as shown in Figure~\ref{fig:disk} and approaches the exact solution numerically as the resolution increases or the grid length decreases.

\subsection{2D calculation with tangential boundary condition}\label{sec:tangentialRect}

To illustrate the tangential boundary condition, we consider a $2\!\times\!2$ square in the center of the $3\!\times\!3$ grid in Figure~\ref{fig:example}. The four solid grid points are the only ones with dual cells that are partially inside the domain. For this particular size, the Hodge stars $S_{0,t}$ and $S_{1,t}$ are both identical matrices, since all the dual cells/edges involved are fully inside the domain. With the following coboundary matrix,
\begin{equation}\label{eq:rectD0t}
D_{0,t} =
\begin{bmatrix}
-1 &  1 &  0 &  0  \\
 0 & -1 &  1 &  0  \\
 0 &  0 & -1 &  1  \\
 1 &  0 &  0 & -1
\end{bmatrix},
\end{equation}
the Hodge Laplacian is identical to the BIG Laplacian, 
\[L_{0,t}^H = L_{0,t}^B =
\begin{bmatrix}
    2 & -1 &  0 & -1 \\
   -1 &  2 & -1 &  0 \\
    0 & -1 &  2 & -1 \\
   -1 &  0 & -1 &  2
\end{bmatrix}.\]
Thus,
\[ \lambda_{0,t}^H= \lambda_{0,t}^B = \{0, 2, 2, 4\}.\]
Note that the exact spectrum below for the Hodge Laplacian contains $0$ with the multiplicity of $\beta_0=1$,
\[\lambda_\text{exact}=\{0, 2.4674, 2.4674, 4.9348\}.\]
When a $1\!\times\!1$ square is given, the BIG Laplacian remains the same as it ignores the change in the Hodge star. The four dual cell size estimate is changed to $0.5^2/2$ whereas the dual edge size estimate is scaled to $0.5$, leading to a doubled frequency. While both spectra are far from the exact spectra, the Hodge star-based calculation is much closer even for such an extremely coarse resolution. It is a simple exercise to see that the discrete $L_{0,t}$ has the exact same spectrum as the discrete $L_{2,n}$ on the dual grid.

If a diagonal edge is added, so to triangulate the square as in the Lagrangian method, we add a row to the matrix from \eqref{eq:rectD0t},
 \begin{equation}\label{eq:rectD0t_triangle}
D_{0,t} =
\begin{bmatrix}
-1 &  1 &  0 &  0  \\
 0 & -1 &  1 &  0  \\
 0 &  0 & -1 &  1  \\
 1 &  0 &  0 & -1 \\
-1 &  0 &  1 &  0
\end{bmatrix}.
\end{equation}
For the diagonal Hodge star, we use cotangent weights of edges. For the new edge, $S_{1,t}$ is $2\cot(\pi/2) = 0$ and all other edges are $\cot(\pi/4) = 1$. While $L_{0,t}^H$ remains the same as above, now
\[L_{0,t}^B =
\begin{bmatrix}
    3  & -1     &   -1  & -1 \\
   -1   & 2  & -1 &        0 \\
    -1  & -1   & 3  & -1 \\
   -1 &        0  & -1  &  2
\end{bmatrix}.\]
Whose third eigenvalue is now considerably skewed,
\[ \lambda_{0,t}^B = \{0, 2, 4, 4\},\]
while the eigenvalues of $L_{0,t}^H$ remain somewhat close to the expected.

\subsection{Impact of Hodge star on combinatorial Laplacian}

\begin{figure}[H]\label{fig:samplings}
    \begin{minipage}[]{.75\linewidth}
    \centering
    \includegraphics[width=\linewidth]{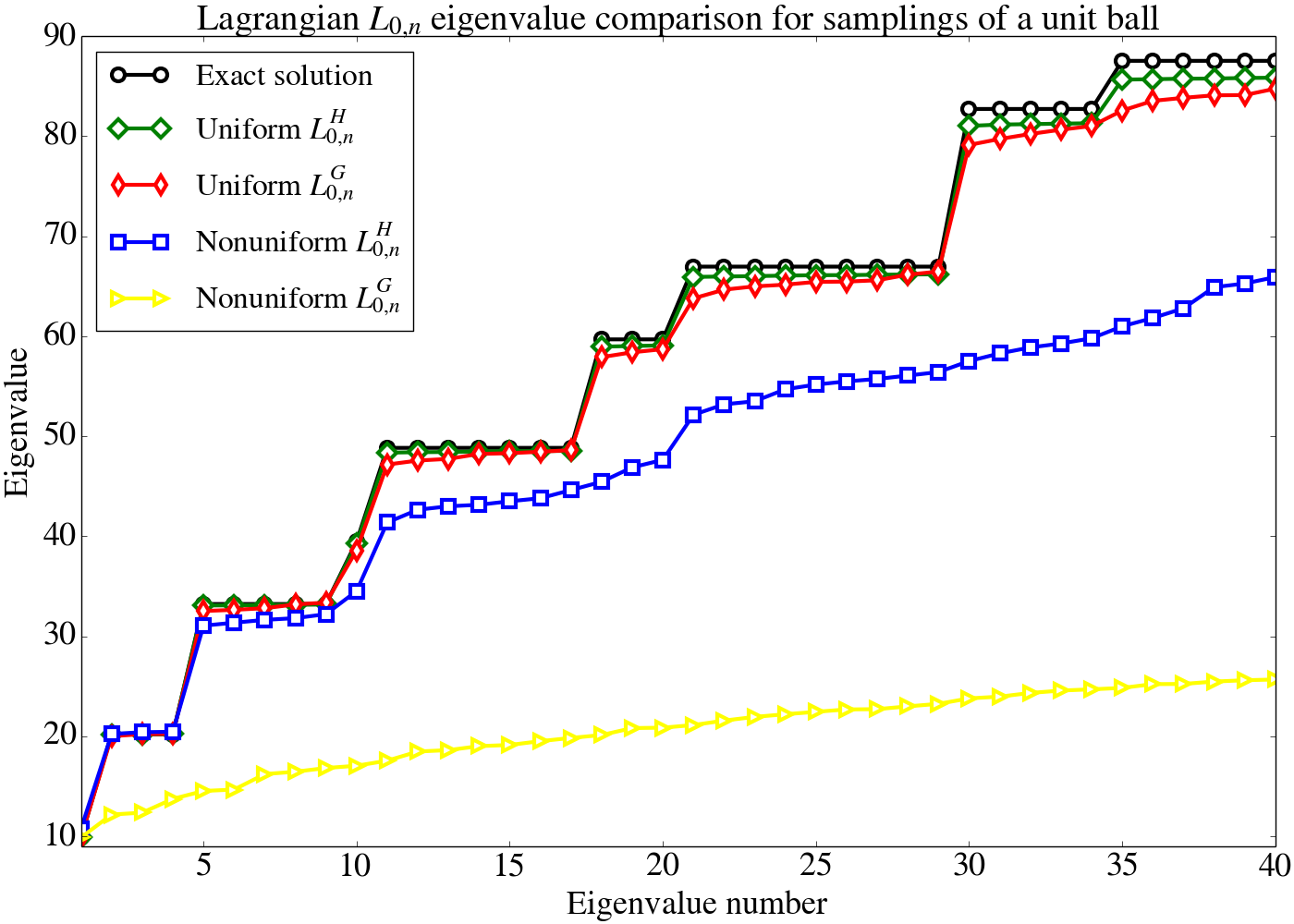}
    \end{minipage}
    \begin{minipage}[]{.24\linewidth}
    \centering
    \includegraphics[width=.97\linewidth]{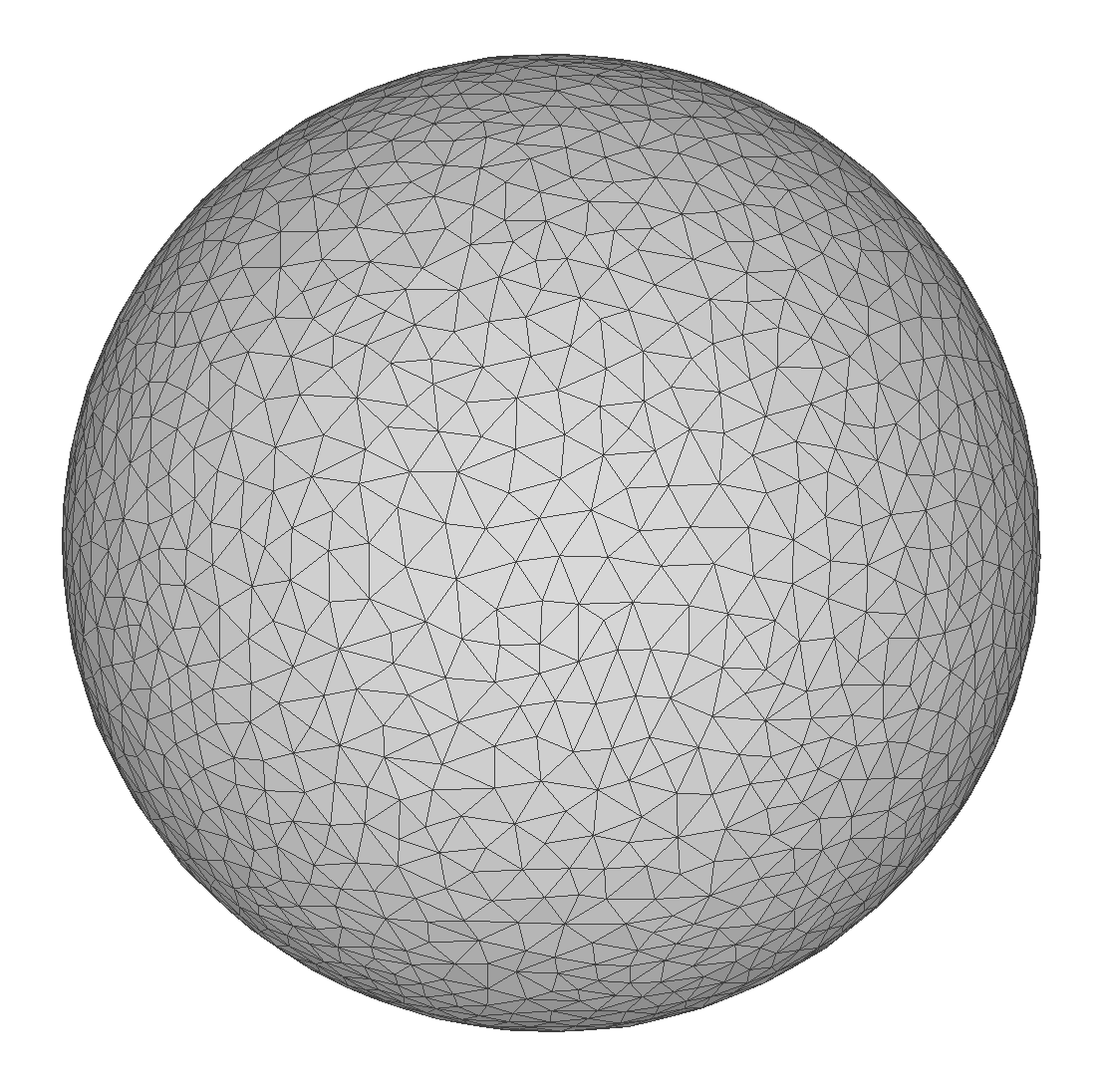} \\
    \includegraphics[width=\linewidth]{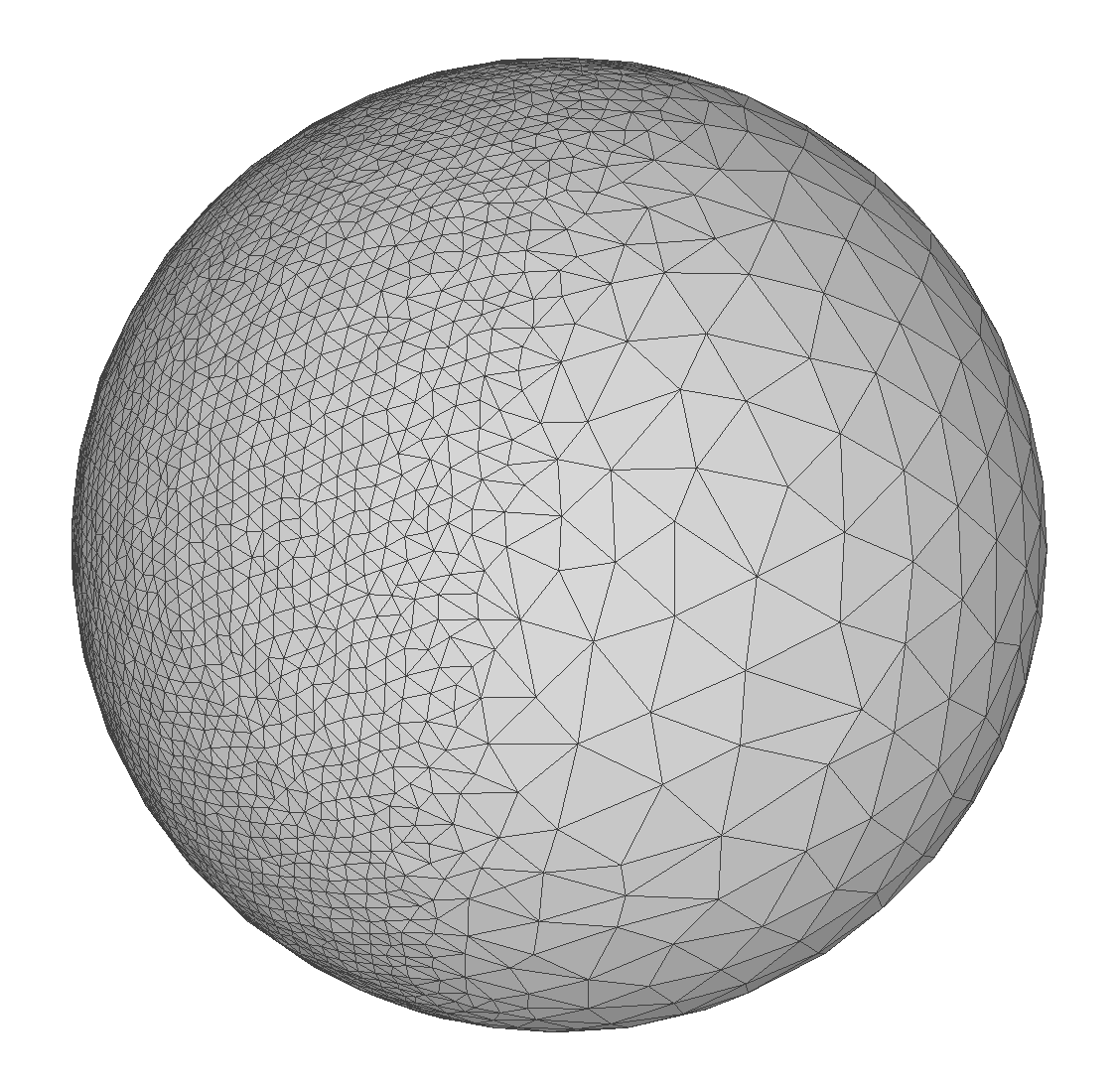}
    \end{minipage}
    \vskip-1em
    \caption{The first 40 eigenvalues of $L_{0,n}$ for the unit ball given uniform and nonuniform samplings. The Lagrangian combinatorial and Hodge Laplacian eigenvalues are shown (left) for boundary surface average edge length approximately 0.07 (top right) and 0.04 (bottom right).}
\end{figure}

Next, we demonstrate the impact of the Hodge star on Lagrangian meshes, i.e., irregularly-sampled point clouds with level set induced connectivity and boundary. When the true Hodge star operator is not a rescaled identity, even for the interior DoF, and not just for the boundary, the combinatorial Laplacian would be nowhere near the expected spectra, as shown in Figure~\ref{fig:samplings}. Whereas, the Hodge star compensates for the irregularity, at least for the first 10 or so eigenvalues.

\vskip-1em
\begin{figure}[H]\label{fig:NTC}
    \begin{center}
    \includegraphics[width=.95\textwidth]{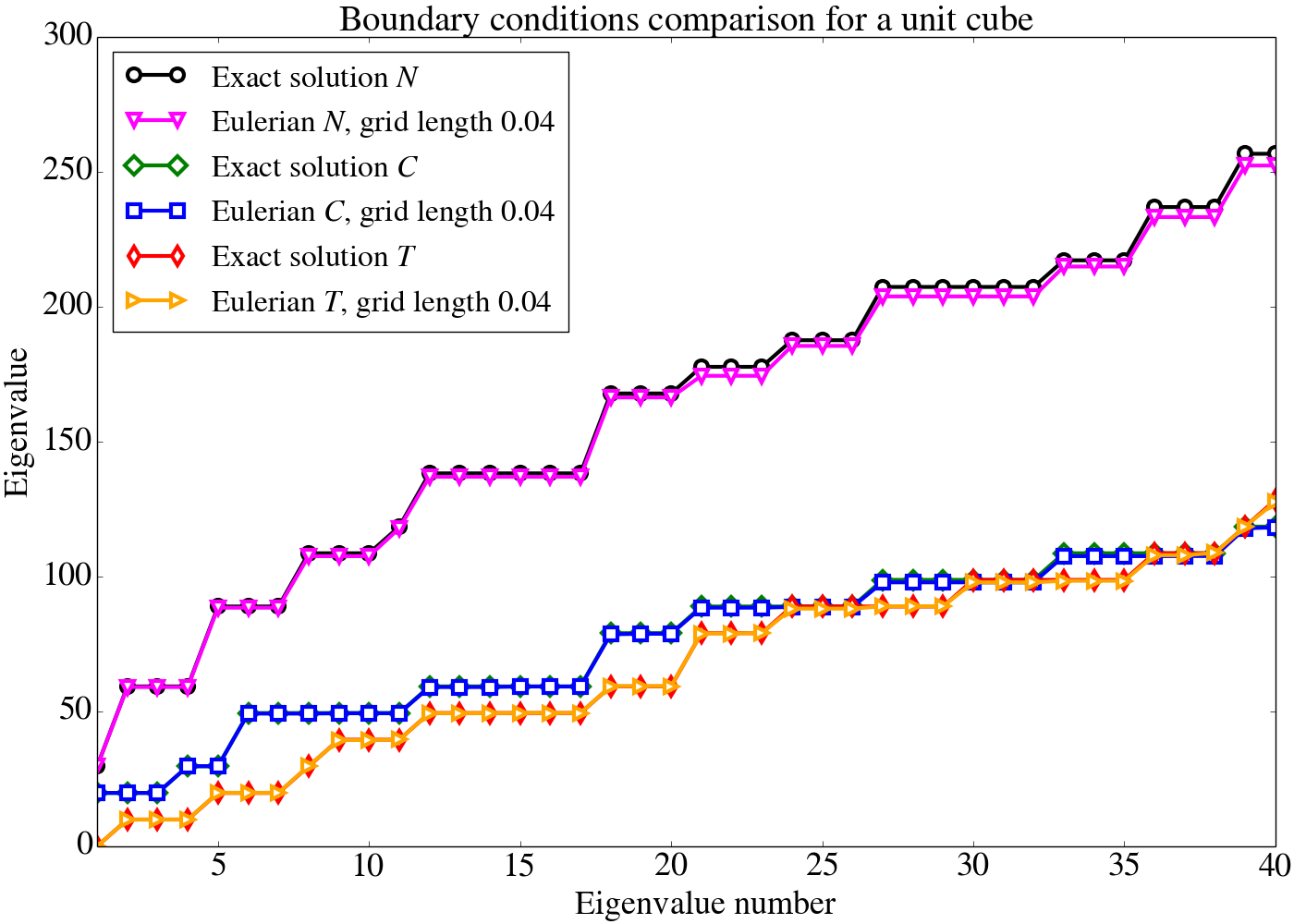}
    \vskip-1em
    \caption{The first 40 eigenvalues of the three spectral sets $N$, $T$, and $C$ for the unit cube. Both exact solutions and Eulerian Hodge Laplacian solutions for grid length 0.04 are shown.}
    \end{center}
\end{figure}

\subsection{Justification for boundary conditions}

Another test where the combinatorial Laplacian fails is for boundary treatment. If the boundary conditions are not correctly handled for the scalar field Laplacian, the variation between $N$ and $T$ can be large as shown in Figure~\ref{fig:NTC}. For the vector Laplacian, the spectra can be either $C\cup N$ or $C\cup T$ depending on the boundary condition, which again are distinct.

\section{Similarities among Laplacians}\label{sec:sim}

Here we demonstrate the similarities of the Eulerian BIG and Hodge Laplacians as well as to that of the Lagrangian representation of the Hodge Laplacian from \cite{zhao2020rham} by comparing their first 40 eigenvalues for $L_{0,n}$, $L_{1,n}$, and $L_{3,n}$. We also include the exact eigenvalues in the sense of the Hodge Laplacian when available and as outlined in the supplementary material.

\subsection{Planar Result}

\begin{figure}[H]\label{fig:disk}
\begin{center}
\includegraphics[width=.95\textwidth]{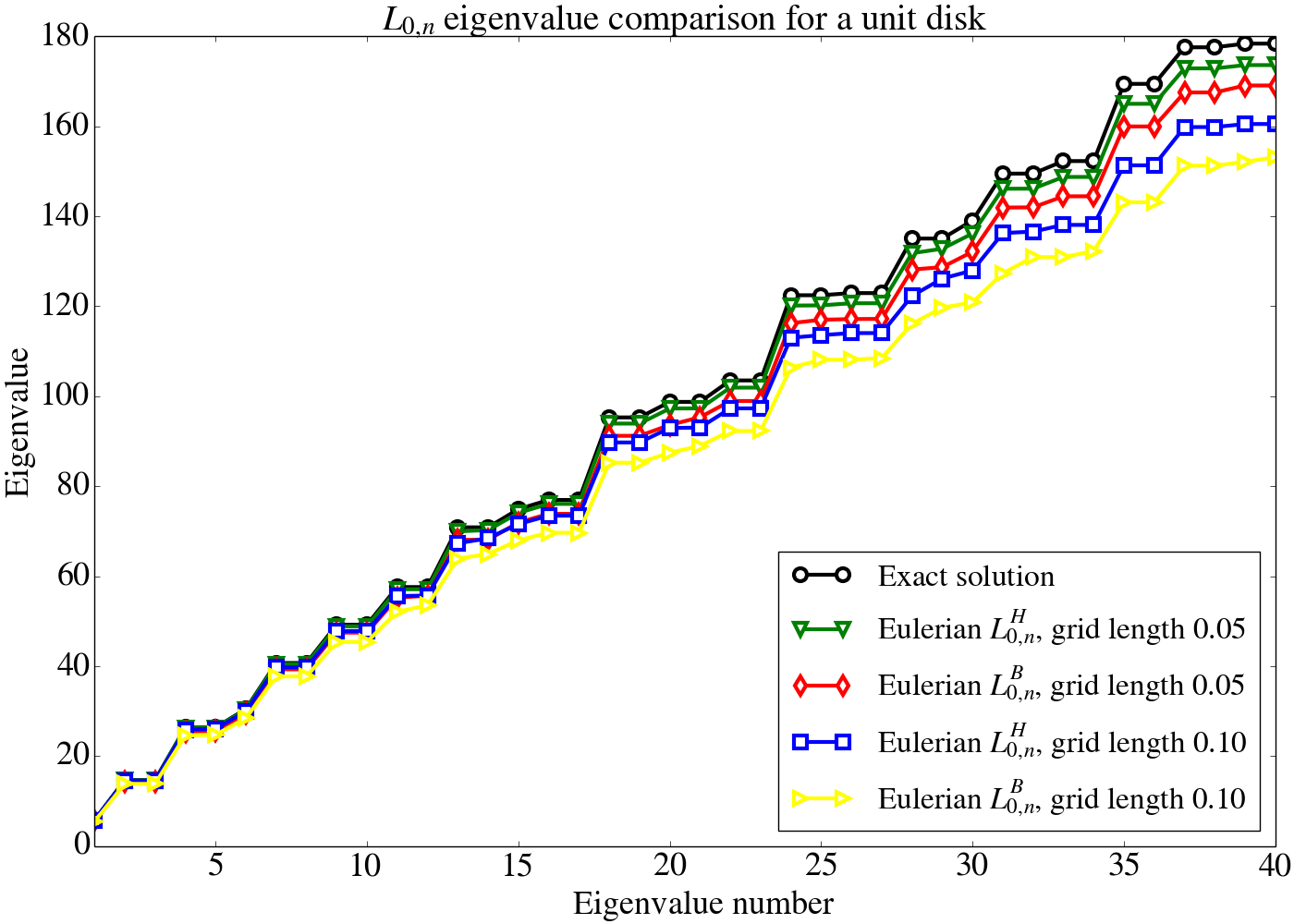}
\vskip-1em
\caption{The first 40 eigenvalues of $L_{0,n}$ for the unit disk. The Eulerian BIG and Hodge Laplacians are shown for two different grid lengths, 0.05 and 0.1.}
\end{center}
\end{figure}

In this 2D example, we computed the 0-Laplacians under normal boundary conditions, i.e., scalar Laplacians with Dirichlet boundary conditions. Figure~\ref{fig:disk} shows the first 40 eigenvalues for grid edge lengths 0.05 and 0.1 for BIG and Hodge Eulerian Laplacians alongside the solutions to the exact eigenvalue problems. The results for the disk in Figure~\ref{fig:disk} show how the addition of the Hodge star, in the case of the Hodge Laplacian, improves the results when compared to the exact solution. In fact, the Hodge Laplacian with a twice as large grid length performs nearly as well as the BIG Laplacian. The exact solutions for the disk are found using Bessel functions. Each eigenvalue $\lambda_{n,k}$ is given by the $k$-th zero of the $n$-th Bessel function and then they are sorted.

\vskip2em

\subsection{Similarities in 3D}\label{3dresults}

Given the duality, in these tests, we compute the 0-, 1-, and 3-Laplacians with normal boundary conditions as their spectra cover the reduced spectral sets. The models we test are the solid cube, ball, torus, and a spherical shell (a ball with a round cavity inside) with one model per each Laplacian. Additional figures are in the supplemental materials. We chose these shapes because they provide examples with trivial and nontrivial homology groups, and their exact solutions are available except for the torus. The Eulerian method performs competitively as opposed to the Lagrangian meshed-based evaluation. The BIG Laplacian, which we extended to handle both boundary conditions, differs from the corresponding Hodge Laplacian only near the boundary in the Cartesian grid. However, our computations demonstrate that the thin layer of difference introduces a great accuracy gain for the Hodge Laplacian.

\begin{figure}[H]\label{fig:models}
    \begin{center}
    \includegraphics[width=.21\textwidth]{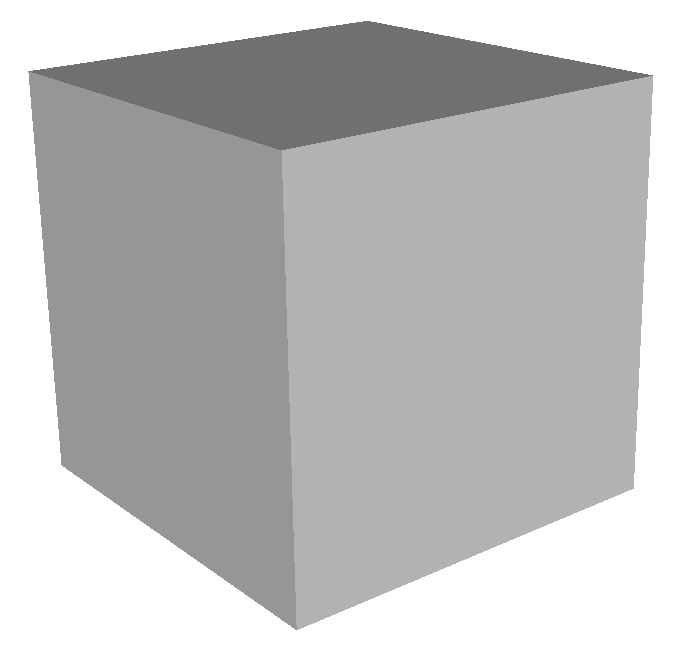}\quad
    \includegraphics[width=.21\textwidth]{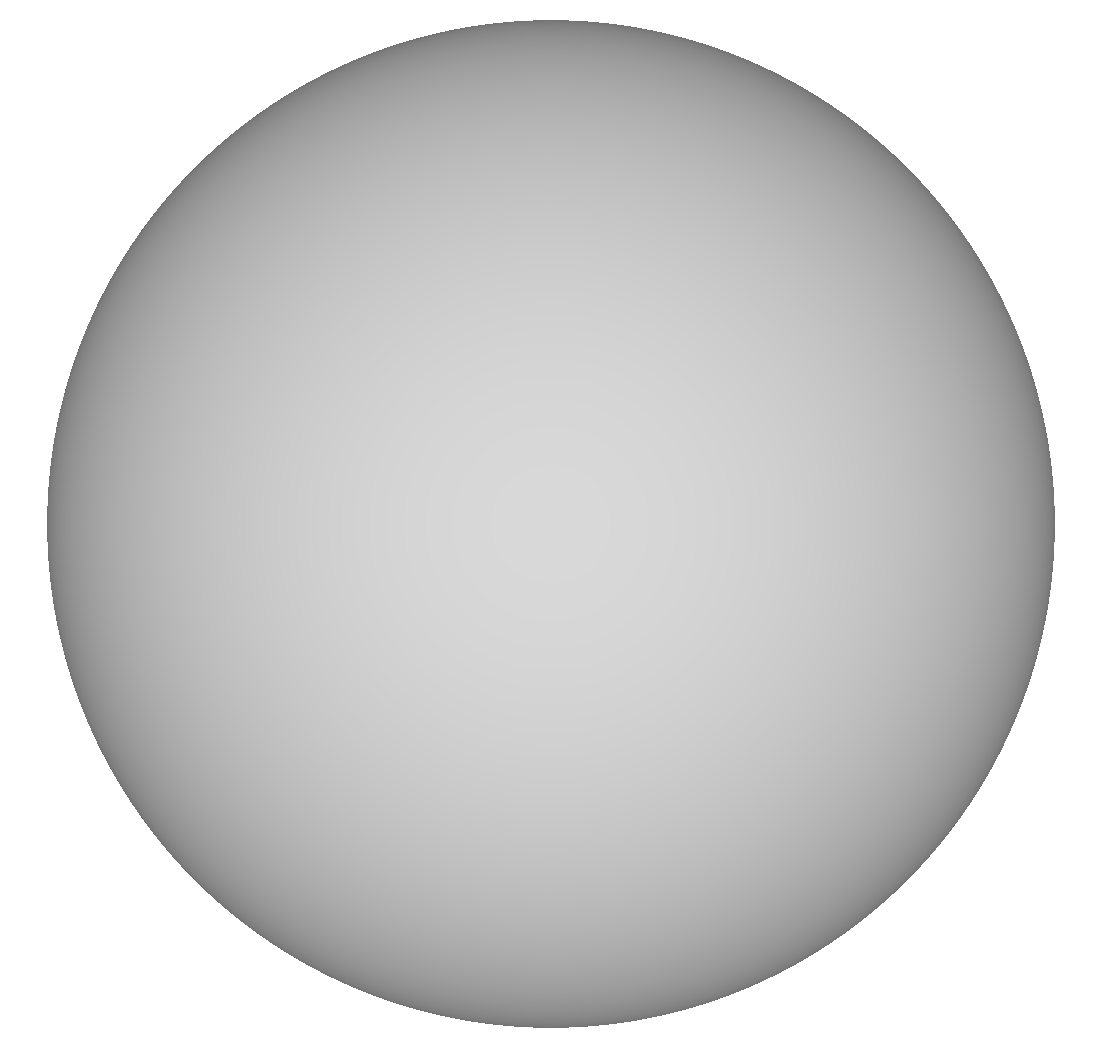}\quad
    \includegraphics[width=.17\textwidth]{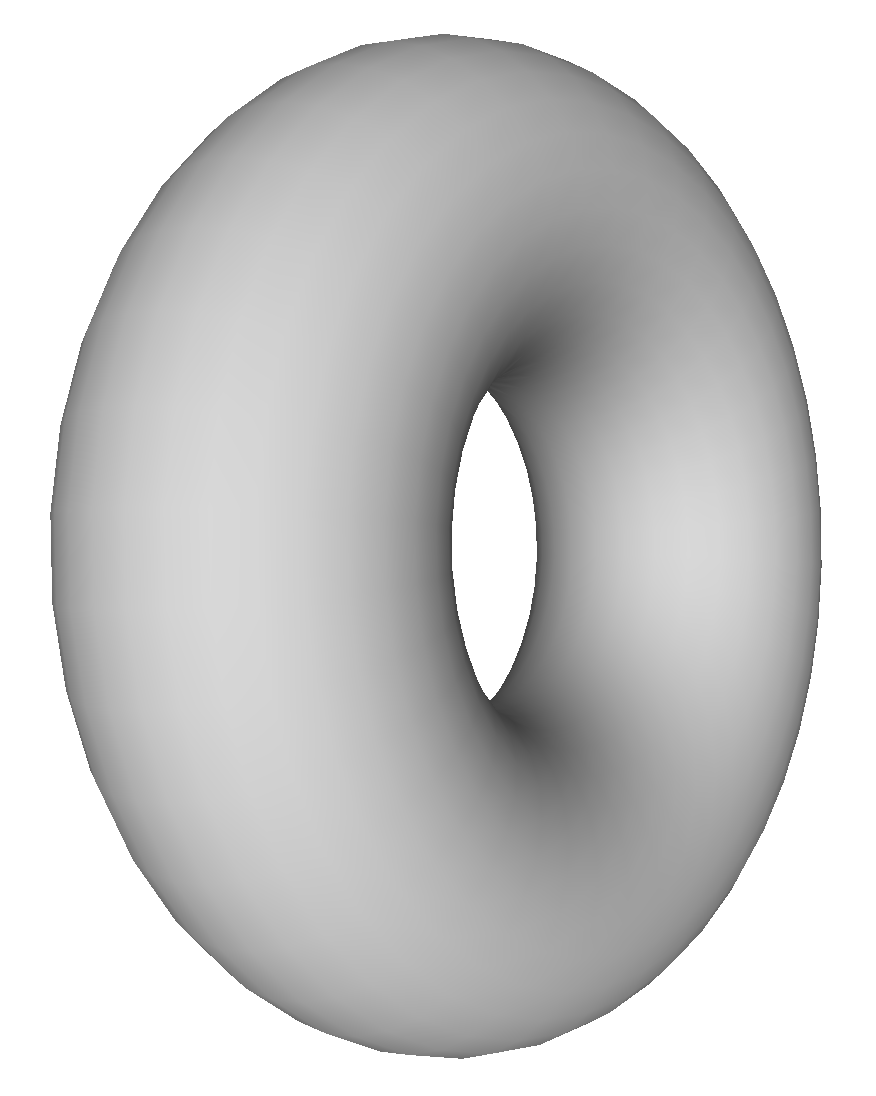}\quad
    \includegraphics[width=.19\textwidth]{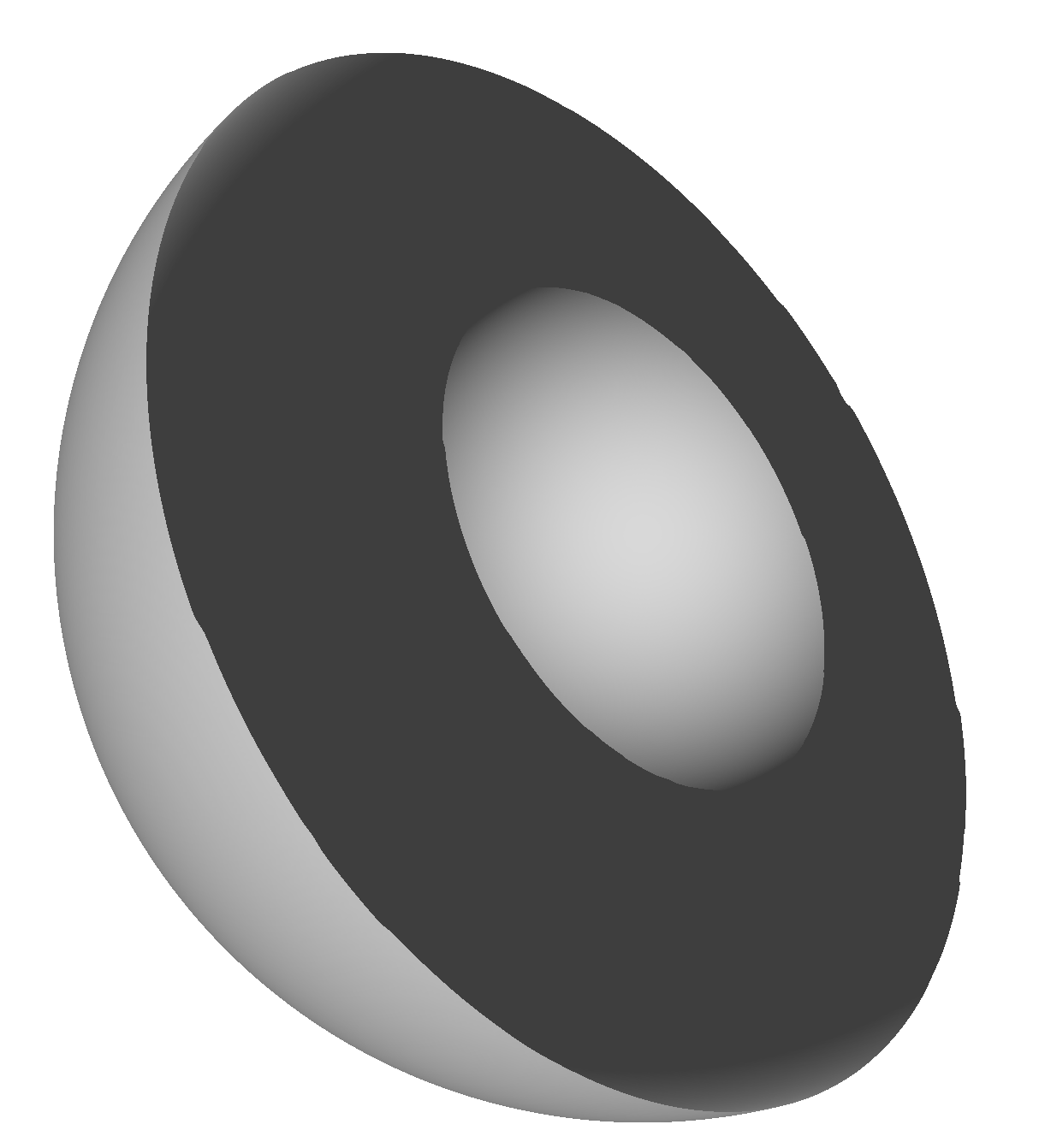}
    \caption{3D models used in Section \ref{3dresults} and the supplementary material.
		 From left to right a cube, ball, torus, and a cut open spherical shell.}
    \end{center}
\end{figure}

\subsubsection{Scalar Laplacian under Dirichlet Boundary Condition $L_{0,n}$}

We compared the spectra of three different discrete Laplacian matrices for scalar fields that vanish on the domain boundary, the Cartesian grid Hodge Laplacian, the BIG Laplacian modified for this boundary condition, and the Hodge Laplacian for a Lagrangian domain discretization, i.e., tetrahedral meshes assuming finite linear elements from \cite{zhao2020rham}. The exact spectra of the continuous Laplacian operators for the cube, the ball, and the spherical shell are given in the supplementary material. For comparisons, we rescale the spectrum from the BIG Laplacian by $1/l_g^2$ as the BIG Laplacian is unitless, but the actual Laplacian carries a unit of one over length squared.

The kernel size of this Laplacian is $0$ as it corresponds to the 0th relative homology with respect to the boundary. Equivalently, it is isomorphic to the kernel of $L_{3,t}$, which has a dimension equaling the Betti number $\beta_3=0.$ Another way to interpret this is that the Laplace equation $\Delta f=0$ on the domain $M$ has a unique solution $f=0$ under the Dirichlet boundary condition $f\vert_{\partial M}=0.$

\begin{figure}[H]\label{fig:cubeL0n}
    \begin{center}
    \includegraphics[width=.95\textwidth]{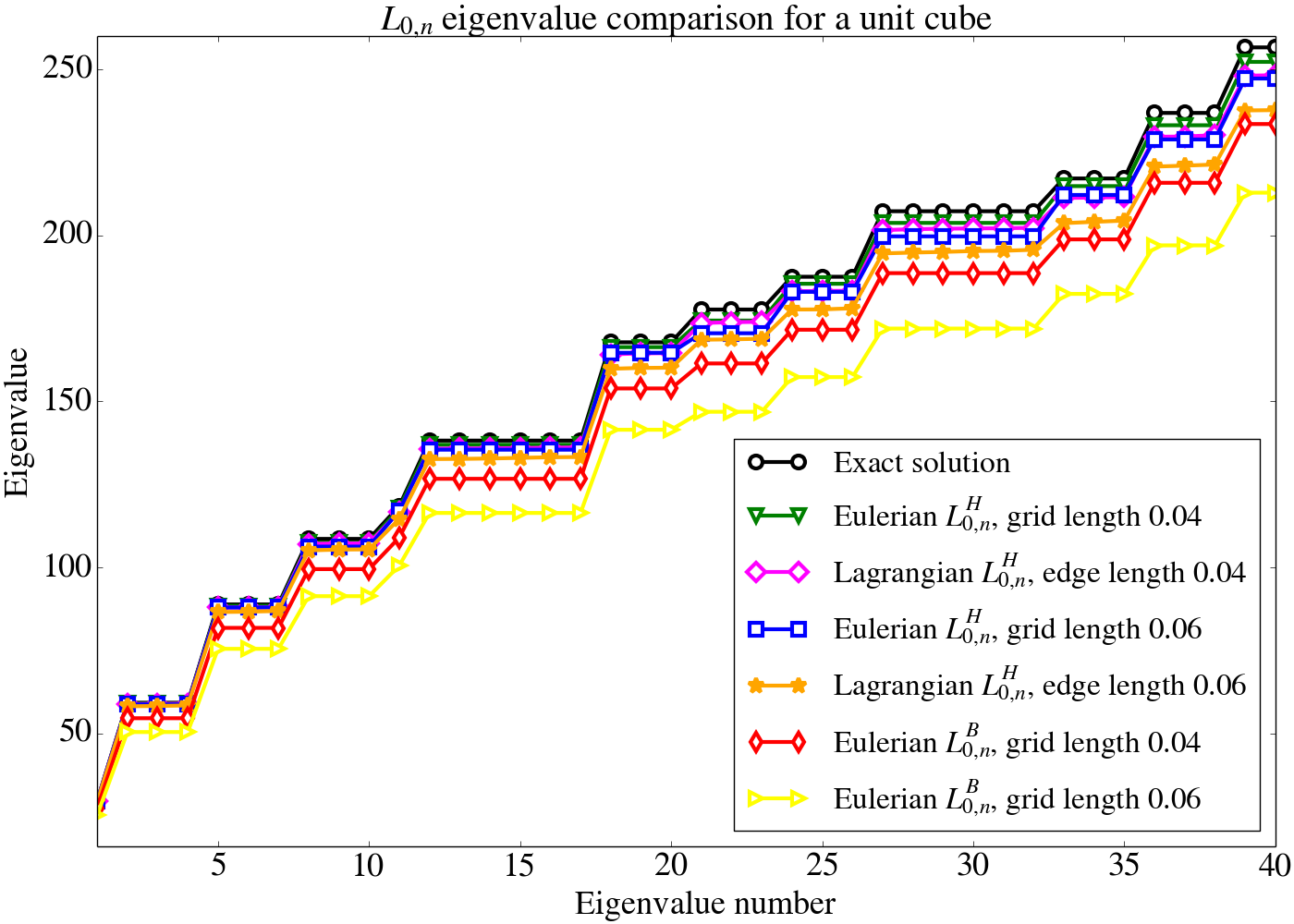}
    \vskip-1em
    \caption{The first 40 eigenvalues of $L_{0,n}$ for the unit cube. The exact solution, Eulerian BIG, and Hodge Laplacians are shown as well as the Lagrangian Hodge Laplacians from \cite{zhao2020rham} for two different grid/edge lengths, 0.06 and 0.04. }
    \end{center}
\end{figure}

For the unit cube in Figure~\ref{fig:cubeL0n}, we see that the Eulerian Hodge Laplacian is the closest to the exact solution, whereas the  BIG Laplacian is the furthest away. We speculate that the better accuracy of the Eulerian Hodge Laplacian compared to the Lagrangian Hodge Laplacian is due to the better alignment of the Cartesian grid with this particular shape. Since the level set function stored near the planar boundary of the cube is parallel to the grid edges, the partial edge lengths interior to the model can be accurately estimated, leading to improved performance of the diagonal Hodge star.
Note that the Lagrangian mesh with a similar resolution requires more DoFs for the scalar field and the Lagrangian Laplacian operator is a denser matrix. The BIG and Hodge Laplacians share the same sparse structure on the grid and are localized to the thin boundary region. However, this small correction to the Hodge star leads to a substantial gain in accuracy.

\vskip-1em

\begin{figure}[H]\label{fig:torusL3n}
\begin{center}
\includegraphics[width=.95\textwidth]{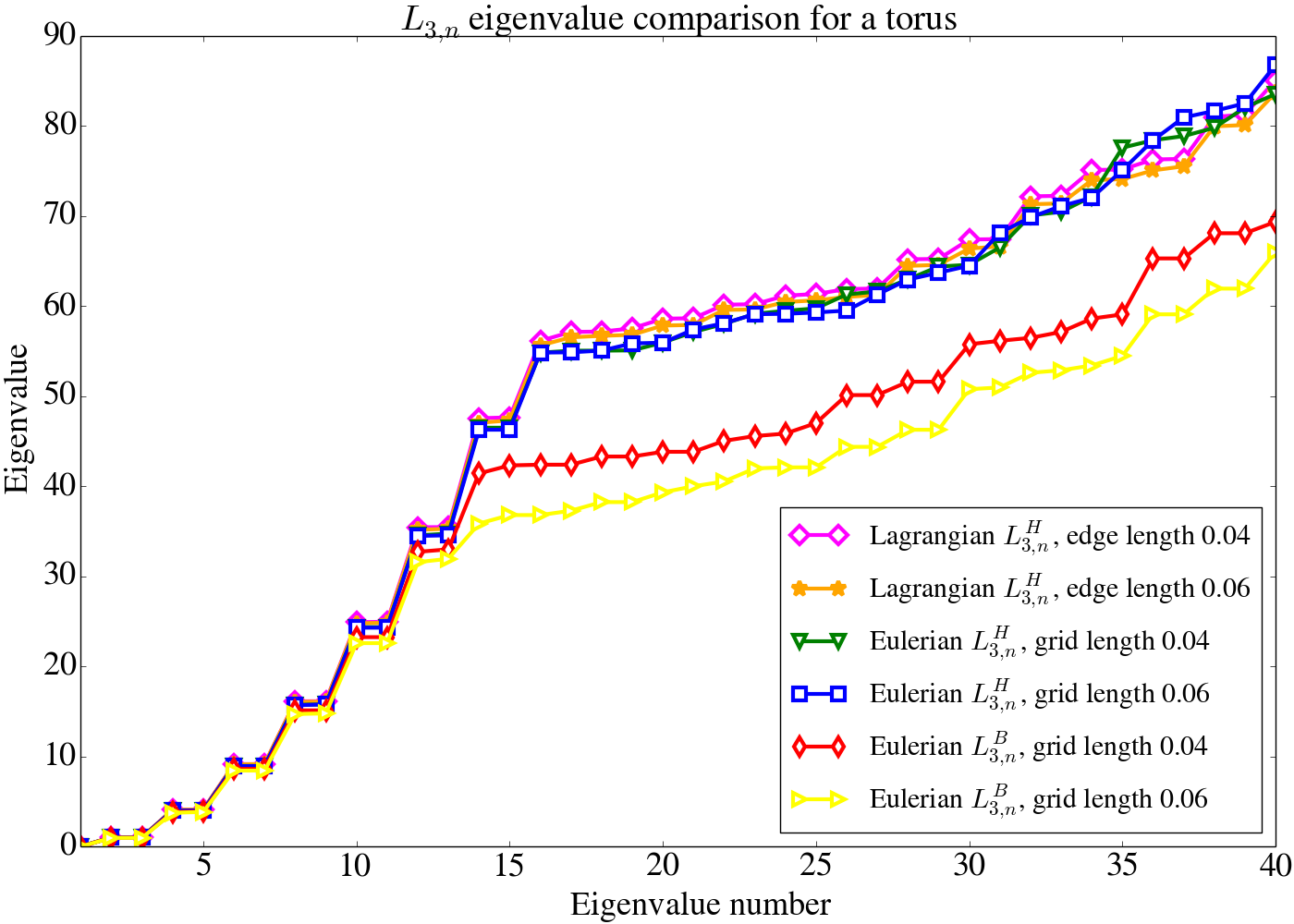}
\vskip-1em
\caption{The first 40 eigenvalues of $L_{3,n}$ for a torus. The Eulerian BIG and Hodge Laplacians are shown as well as the Lagrangian Hodge Laplacians from \cite{zhao2020rham} for two different grid/edge lengths, 0.06 and 0.04. No exact solution was computed for the torus.}
\end{center}
\end{figure}

\subsubsection{Scalar Laplacian under Neumann Boundary Condition $L_{3,n}$}

We compared the spectra of two different discrete Laplacian matrices for scalar fields that have a $0$ normal derivative on the domain boundary.
The exact spectra of the continuous Laplacian operator for the ball and other models are given in the supplementary materials. The kernel size of this Laplacian corresponds to the dimension of the 3rd relative homology with respect to the boundary, or equivalently, the 0th homology, i.e., the number of connected components $\beta_0.$ We can also see that a constant scalar field is the only one satisfying the condition, which is still true with the discrete Hodge or BIG Laplacians. As expected and shown in Figure~\ref{fig:torusL3n}, the graph version consistently produces an accuracy far worse than the lower resolution results of the Hodge versions. Nevertheless, the kernel sizes (topological invariants) are correct. The other models are included in the supplementary material.

\begin{figure}[H]\label{fig:shellL1n}
\begin{center}
\includegraphics[width=.95\textwidth]{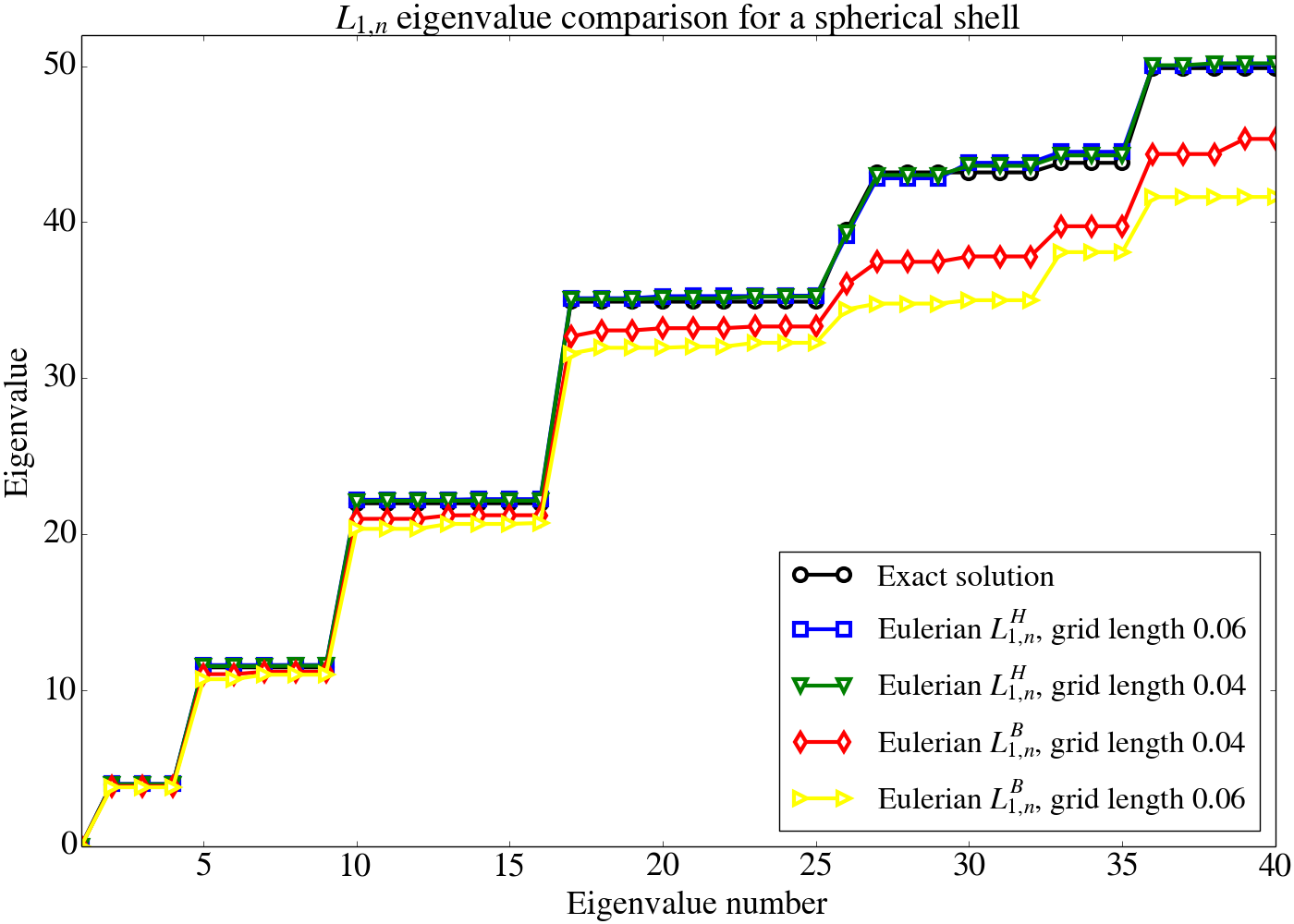}
\vskip-1em
\caption{The first 40 eigenvalues for $L_{1,n}$ for a spherical shell with outer radius 1.0 and inner radius 0.5. The exact solution, Eulerian BIG, and Hodge Laplacians are shown for two different grid  edge lengths, 0.06 and 0.04.}
\end{center}
\end{figure}

\subsubsection{Vector Laplacian under Normal Boundary Condition $L_{1,n}$}

In the test for the Laplacian of vector fields that are normal to the domain boundary (along with the two Robin conditions to make the kernel finite), we found similar behaviors for the Laplacians as shown in Figure~\ref{fig:shellL1n}. The exact spectra of the continuous Laplacian operator for the spherical shell are again given in supplementary materials. The kernel size of this Laplacian corresponds to the dimension of the 1st relative homology with respect to the boundary, or equivalently, the 2nd homology, i.e., the number of cavities $\beta_2$, of which the spherical shell has a nontrivial kernel. Our test confirms that both Hodge and BIG Laplacians do produce the right kernel. While this is the only nontrivial test, it can be explained by Algorithm.~\ref{alg:buildops}. Any $k$-cell is included whenever one of its vertices is inside, so if a dual cell is included, all its vertices, edges, and faces are included. Thus, the topology is the same as that of the voxelization of the original shape on the dual grid, as long as the grid has a sufficiently fine resolution.

\section{Concluding Remarks}

Although both the combinatorial Laplacian and Hodge Laplacian can reveal the topological dimension and geometric shape of data, they are defined on discrete point clouds and continuous manifolds, respectively. ``Hodge Laplacian on graphs''  \cite{lim15} emphasizes the apparent structural similarity of combinatorial Laplacian  and Hodge Laplacian in algebraic topology. However, combinatorial Laplacians defined on simplicial complexes are conceptually incompatible with Hodge Laplacians on manifolds with or without boundary. Although in computational mathematics, using discrete vector calculus (DEC),  gradient, curl, and divergence can be defined on discrete grids i.e., point clouds with boundary or regular meshes, their realization on simplicial complexes is not possible because generic simplicial complexes themselves in the combinatorial Laplacian setting do not admit a functional space as in the cellular sheaf theory setting.
Additionally,  it is impossible for combinatorial Laplacians to have the Hodge decomposition of a vector field into physically relevant harmonic, curl-free, and divergence-free components as shown in Figure \ref{fig:decomposition}. Instead, combinatorial Laplacian decompositions reflect the connectivity among simplices of different dimensions. As a result, the two approaches cannot be further compared.

This work introduces Boundary-Induced Graph (BIG) Laplacian to bring discrete Laplacians and continuous Hodge Laplacian on an equal footing for detailed analysis and comparison. BIG Laplacians of various topological dimensions are defined on Cartesian domains with proper boundary conditions (i.e., Dirichlet and Neumann) to deliver the correct topological dimensions of data.

In the context of computing Betti numbers and the smallest eigenvalues, the BIG Laplacian may perform favorably and is simpler to compute, with proper modifications for boundary conditions. As shown, given uniform sampling, either a regular grid or uniform tetrahedral mesh, the spectra of BIG  Laplacians have worse ``accuracy'' compared to the Hodge Laplacian spectra. However, if the input is irregularly sampled, which could be the case with real-world data, the Hodge star would be indispensable to account for nonuniform geometric quantities. The Hodge Laplacian preserves higher frequency spectral information, such as for shape and spectrum learning tasks. Furthermore, while the duality between normal and tangential boundary conditions indicates that the normal boundary condition of $k$-forms can be handled with the clique-based combinatorial Laplacian for $(n\!-\!k)$-forms, leading to correct kernel dimensions, only the correct boundary treatment through primal and dual grids can handle mixed boundary conditions, which are common in practical problems.

In practical applications, BIG Laplacians stand out as an independent formulation for the analysis of volumetric data, such as those from Cryo-Electron Microscopy (Cryo-EM),  Magnetic Resonance Imaging (MRI), computer vision, and the solutions of Partial Differential Equations (PDEs). Its generalization to evolving manifolds, namely persistent BIG Laplacians (PBLs), in analogous to evolutionary de Rham-Hodge theory \cite{chenevol}, will have much potential for machine learning modeling and prediction.

 \section*{Code availability}
BIG Laplacian code can be found at \\ \begin{center}\href{https://github.com/eribandogros/BIGLaplacians}{https://github.com/eribandogros/BIGLaplacians}.\end{center}

\section*{Supplementary material}
The Supplementary material is available for
\begin{enumerate}
	\item[S1]   Exact Laplacian spectra of elementary shapes\\
	S1.1. Spectra of Laplacians on the unit ball\\
	S1.2. Vector field Laplacian in spherical polar coordinates\\
	S1.3. Spectrum of vector field Laplacian on ball\\
	S1.4. Spectra on spherical shells\\
	S1.5. Spectra on cuboids
	
	\item[S2]   Additional Results\\
	S2.1. Planar Result\\
	S2.2. 3D Results\\
	S2.2.1. Scalar Laplacian under Dirichlet Boundary condition $L_{0,n}$\\
	S2.2.2. Scalar Laplacian under Neumann Boundary condition $L_{3,n}$\\
	S2.2.3. Vector Laplacian under Normal Boundary condition $L_{1,n}$
\end{enumerate}

\section*{Acknowledgment}
This work was supported in part by  NSF grant IIS-1900473 and  NIH grants  GM126189 and R01AI164266.

 \bibliographystyle{siamplain}
\bibliography{refs_new}
\end{document}